\newcommand{\lapack}{LAPACK }
\newcommand{\blas}{BLAS }
\newcommand{\scalapack}{ScaLAPACK }
\newcommand{\MBF}[1]{{#1}}

\newcommand{\BW}{k}

\newcommand{\RECLEVELS}{r}
\newcommand{\udots}{\reflectbox{$\ddots$}}
\newcommand{\matmid}[1]{\multirow{2}{*}{$#1$}  \\ \hspace{0cm}
} 
\newcommand{\matmidt}[1]{\multirow{2}{*}{$#1$} \\ \hspace{0cm} } 
\newcommand{\T}[1]{\rule{0pt}{#1ex}}

\documentclass[acmtoms]{acmsmall} % Aptara syntax

\usepackage{arydshln}

\usepackage[ruled]{algorithm2e}

\usepackage{listings}
\usepackage{setspace}
\usepackage{mdframed}
\SetAlFnt{\small}
\SetAlCapFnt{\small}
\SetAlCapNameFnt{\small}
\SetAlCapHSkip{0pt}
\IncMargin{-\parindent}

\usepackage{graphicx}
\usepackage{tikz}
\usepackage{amsmath}
\usepackage{multicol}
\usepackage{multirow}
\usepackage[thinlines]{easybmat}
\usepackage{caption}

\definecolor{darkred}{RGB}{136,28,28}
\definecolor{white}{RGB}{255,255,255}
\definecolor{blueUMASS}{RGB}{8,52,79}
\definecolor{lightblueUMASS}{RGB}{149,183,215}
\definecolor{redUMASS}{RGB}{136,28,28}
\definecolor{yellowUMASS}{RGB}{251,214,104}
\definecolor{orangeUMASS}{RGB}{161,108,15}
\definecolor{greenUMASS}{RGB}{93,127,26}
\definecolor{grayUMASS}{RGB}{160,160,160}
\definecolor{tanUMASS}{RGB}{190,179,149}
\definecolor{greyUMASS}{RGB}{160,160,160}
\makeatletter
\DeclareRobustCommand{\rvdots}{%
  \vbox{
    \baselineskip4\p@\lineskiplimit\z@
    \kern-\p@
    \hbox{.}\hbox{.}\hbox{.}
  }}
\makeatother

\begin{document}
%%%%%%%%%%%%%%%%%%%%%%%%%%%%%%%%%%%%%%%%%%%%%%%%%%%%%%%%%%%%%%%%%%%%%%%%%%%%%%%%%%%
% Page heads
\markboth{Spring, Polizzi, Sameh}{A Feature Complete SPIKE Banded Algorithm and Solver}

% Title portion
\title{A Feature Complete SPIKE Banded Algorithm and Solver}
\author{
BRAEGAN S. SPRING and ERIC POLIZZI
\affil{University of Massachusetts, Amherst}
AHMED H. SAMEH
\affil{Purdue University}
}
% NOTE! Affiliations placed here should be for the institution where the
%       BULK of the research was done. If the author has gone to a new
%       institution, before publication, the (above) affiliation should NOT be changed.
%       The authors 'current' address may be given in the "Author's addresses:" block (below).
%       So for example, Mr. Abdelzaher, the bulk of the research was done at UIUC, and he is
%       currently affiliated with NASA.

\begin{abstract}
  New features and enhancements for the SPIKE banded solver are presented.
  Among all the SPIKE algorithm versions, we focus our attention on
  the recursive SPIKE technique which provides
  the best trade-off between generality and parallel efficiency, but was known for
  its lack of flexibility. Its application  was essentially limited to 
  power of two number of  cores/processors. This limitation is successfully
  addressed in this paper.  In addition, we present a new
  transpose solve option,
  a standard feature of most numerical solver libraries which has
  never been addressed by the SPIKE algorithm so far.
  A pivoting recursive SPIKE strategy is finally presented as an alternative to
  non-pivoting scheme for systems with large condition numbers.
  All these new enhancements participate
  to create a feature complete SPIKE algorithm and a
  new black-box SPIKE-OpenMP package that significantly outperforms the performance
  and scalability obtained with other state-of-the-art banded solvers.
\end{abstract}

%\category{C.2.2}{Category1}{Category2}

%\terms{TheTerms}

%\keywords{SPIKE, banded matrices, linear system solver}

%\acmformat{Authors, date, PaperTitle}

% At a minimum you need to supply the author names, year and a title.
% IMPORTANT:
% Full first names whenever they are known, surname last, followed by a period.
% In the case of two authors, 'and' is placed between them.
% In the case of three or more authors, the serial comma is used, that is, all author names
% except the last one but including the penultimate author's name are followed by a comma,
% and then 'and' is placed before the final author's name.
% If only first and middle initials are known, then each initial
% is followed by a period and they are separated by a space.
% The remaining information (journal title, volume, article number, date, etc.) is 'auto-generated'.

%\begin{bottomstuff}
%This work is supported by NSF grant CCF-\#1510010.
%\end{bottomstuff}

\maketitle

\newpage 
\tableofcontents

\section{Introduction} \label{Motivation}
Linear systems (i.e. find $X$ solution of $AX=F$ for a given square matrix $A$ and right hand side vectors $F$)
are a fundamental tool, frequently used to express our understanding of the natural and engineering world.
Because of the importance of linear systems in applications,
high quality linear algebra software is a cornerstone of computational science.
Two well known examples of software for performing dense and banded linear algebra are \blas (Basic Linear Algebra Subprograms) and \lapack
(Linear Algebra PACKage) \cite{LAPACK}.
These collections of subroutines provide a consistent interface to high performance linear algebra building blocks across hardware platforms and operating systems.

Many recent improvements in available computational power have been driven by increased use of parallelism. 
The development of new parallel algorithms for solving linear systems
aims at achieving scalability and performance over \lapack LU algorithms on either shared memory or distributed memory architectures.
In shared memory systems, the parallelism in \lapack LU can directly benefit from the threaded implementation of the low-level \blas routines.
In order to achieve further scalability improvement, however, it is necessary to move to a higher level of parallelism based on divide-and-conquer techniques. 
The latter are mandatory with the use of distributed memory systems but they are also becoming increasingly important if one aims at fully exploiting 
shared memory machines composed of a large number of cores. The LU factorization paradigm could be adapted to address a high-level parallelism implementation
as it is the case for the algorithms proposed in the \scalapack library package \cite{ScaLAPACK}. However, in many instances, it can become more advantageous to 
design algorithms that are inherently better suited for parallelism such as the SPIKE algorithm for solving banded linear systems.

This paper is focusing on one particular class of sparse linear systems that are banded. For example, a tridiagonal matrix is a particular
banded matrix with bandwidth of size $b=3$. In practice, $b$ could be much larger, and systems are considered banded if $b<<n$ where $n$ is the size of $A$.
The systems could either be dense or sparse within the band, but only the dense banded case is considered here (the band may explicitly include zero elements).
Very often, banded systems arise after a general sparse system is reordered in some fashion \cite{Cuthill:1969} or
they can naturally arise from applications (e.g. \cite{PolizziJCP:2005}). In other instances, they are
constructed as effective preconditioners to general sparse systems where they are solved via iterative methods \cite{Murat:2010}.

SPIKE is a very effective banded solver
  which can 
  significantly outperform  the ScaLAPACK package on distributed memory systems, as well as LAPACK on shared memory systems. A SPIKE-MPI package was released
in collaboration with Intel in 2008 \cite{Polizzi:2006,SPIKEINTEL,Polizzi:2011},
and a SPIKE-OpenMP solver was completed in 2015 and included
into the distribution of FEAST eigenvalue solver v3.0 \cite{PolizziFEAST,FEASTUG,FEAST} (where
SPIKE is used as a kernel for solving banded eigenvalue problems). GPU implementations of SPIKE have also been proposed by other authors \cite{Chang:2012,Li:2014}.

This work  presents essential enhancements to
  the SPIKE algorithm that are required to achieve a feature complete SPIKE library package.
The development of a competitive library package must not only be motivated by
good performance results, it should consider integrating all the main features offered by common packages.
Among the large number of variants available for SPIKE, we are focusing our efforts
to expand the capabilities of the recursive (parallel) SPIKE algorithm. The recursive scheme demonstrates parallel efficiency and 
is applicable to both diagonally and non-diagonally dominant systems. However, it lacked the flexibility to adapt to some key situations.
In this work, new features and usability enhancements for recursive SPIKE will be considered in order to address the issues listed below.
\begin{enumerate}
\item In practice, the standard SPIKE recursive scheme
is prone to potential waste of parallel resources if the number of cores/processors is not a power of two. For instance, if SPIKE runs
on 63 cores, then only 32 would be effectively used (i.e. the lowest nearest power of two). Here, this restriction is fully removed using
a new flexible partitioning scheme and load balancing strategy that will be presented in Section~3.
\item Most library solvers include the 'transpose solve' option as a standard feature. The same factorization of the matrix $A$ can then be used
  to solve either $AX=F$ or $A^TX=F$ (i.e. there is no need to factorize $A^T$). This feature is important in many practical situations including
  the efficient use of non-transpose free iterative solvers (where $A$ is a preconditioner), and the capability to achieve a $\times 2$ speedup while solving
  complex Hermitian and non-Hermitian eigenvalue problems using FEAST \cite{Kestyn:2016}. The transpose solve option for the SPIKE algorithm is successfully derived 
  in Section~4.
\item The SPIKE recursive scheme is usually associated with a non-pivoting factorization strategy applied to each matrix partition.
The non-pivoting option in SPIKE helps maintaining the
banded structure of the matrix, which simplifies the implementation of the algorithm
and improves performance of the factorization stage. For systems with large condition numbers, however,
partial pivoting may become a necessity for obtaining low residual solutions (without the need to perform iterative refinements).
An efficient pivoting scheme for the recursive SPIKE is presented in Section~5.
\end{enumerate}

All these new enhancements participate to create a feature complete SPIKE algorithm which can be utilized to implement a SPIKE-MPI or SPIKE-OpenMP
library. Without loss of generality (since both MPI and OpenMP SPIKE are concerned), the presentation terminology
and all numerical results are considering a SPIKE OpenMP implementation and the use of threading.
A broader impact of this work has been the development and released of a new stand-alone
SPIKE-OpenMP package (v1.0)  \cite{spike}.
To the extent possible, this solver has been designed as an easy to use, `black-box' replacement to the standard \lapack banded solver.
For example, the library includes support for single and double precision arithmetic using real or complex system matrices.  
Sections 4 to 6 of this paper are accompanied with extensive numerical experiments that demonstrate that the SPIKE solver
significantly outperforms the performance and parallel scalability obtained using the \lapack banded solver in Intel-MKL.
The basic SPIKE algorithm using the recursive scheme is first summarized in Section~2.

\section{SPIKE background} 

The SPIKE algorithm is a domain decomposition method for solving block tridiagonal matrices. 
It can be traced back to work done by A. Sameh and associates on block tridiagonal system in the late seventies \cite{Chen:1978,Sameh:1978,Gallivan2012}.
The central idea in SPIKE departs from the traditional $LU$ factorization 
  with the introduction a new $\MBF{DS}$ factorization which is better suited for parallel implementation as it naturally leads to lower
  communication cost. Several enhancements and variants of the SPIKE algorithm have since been proposed by Sameh and coauthors
  in \cite{Dongarra:1984,Lawrie:1984,berry:1988,Sarin:1999,Polizzi:2006,Polizzi:2007,Murat:2009,Maxim:2010,Murat:2010,Murat:2011}.
Parallelism is extracted by decoupling the relatively large blocks along the diagonal, solving them independently,
and then reconstructing the system via the use of smaller reduced systems. 
There are a number of versions of the SPIKE algorithm, which handle the specifics of those steps in different ways. 
Two main families of algorithms have been proposed in recent years \cite{Polizzi:2006,Mikkelsen:2009,Mendiratta:2011}:
(i) the truncated SPIKE algorithm for diagonally dominant systems; and (ii) the recursive SPIKE algorithm for general non-diagonally dominant systems.
This paper describes improvements to the recursive SPIKE algorithm for solving banded matrices which can either be diagonally or non-diagonally dominant.

\subsection{Central concept of SPIKE} \label{basics_section}
This section presents the basic SPIKE algorithm that will be used to build upon. 
The goal is to find $\MBF{X}$ in the equation
\begin{equation}
{A}{X}={F}, 
\end{equation}
where $\MBF{A}$ is a banded, $n\times n$ matrix.
For clarity, the number of super and sub-diagonals is assumed to be the same and equal to $\BW$. The matrix bandwidth is $b=2\BW+1$ where $\BW$ denotes then
  the ``half-bandwidth''.
The modifications to allow for matrices with non-symmetric bandwidth consist primarily of padding various small submatrices in the SPIKE reduced system with zeroes.
The size of matrices $\MBF{F}$ and $\MBF{X}$ is $n \times n_{rhs}$.

The banded structure may be exploited to enable a domain decomposition. 
$\MBF{A}$ is partitioned along the diagonal into $p$ main diagonal submatrices $A_i$ and their interfaces, as follows:

\begin{equation} \label{partitioning}
\MBF{A} = 
\left[
\begin{array}{ccccc}
\MBF{{A_1}}&\MBF{{B_1}}&&& \\
 \MBF{C_2}&\MBF{{A_2}}&\MBF{B_2}&& \\
 &\ddots&\ddots&\ddots&\ddots \\
 &&&\MBF{C_p}&\MBF{{A_p}}
\end{array}
\right].
\end{equation}

Each $\MBF{A_i}$ is a square matrix of size $n_i$. 
Because the matrix is banded, $\MBF{B_i}$ and $\MBF{C_i}$ can be considered tall and narrow matrices  of size $n\times k$ which contain primarily zeroes
i.e.

\begin{equation}
\label{BandChats}
\MBF{B_i} = 
\left[ \begin{BMAT}(b,18pt,18pt){c}{bc:c}
\matmid{\MBF{0}}\\ 
 \MBF{\hat{B}_i}
  \end{BMAT} \right];
\quad
\MBF{C_i} =
\left[ \begin{BMAT}(b,18pt,18pt){c}{b:cc}
     \MBF{\hat{C}_i} \\
\matmid{\MBF{0}}
  \end{BMAT} \right],
\end{equation}
where $\MBF{\hat{B_i}}$ and $\MBF{\hat{C_i}}$ are small dense  square matrices of size $k$.

We can now factorize the $\MBF{A}$ matrix into the $\MBF{D}$ and $\MBF{S}$ matrices.
$\MBF{D}$ contains the diagonal blocks of the matrix $\MBF{A}$.
$\MBF{S}$ (a.k.a. the spike matrix) relates the partitions to one another as follows:
\begin{equation}
\MBF{A} = 
\MBF{DS} = 
\left[\begin{array}{ccccc}
\MBF{D_1}&&&& \\
 &\MBF{D_2}&&& \\
 &&\ddots&& \\
 &&&&\MBF{D_p}
\end{array}\right]
\left[\begin{array}{ccccc}
\MBF{I_1}&\MBF{V_1}&&& \\
 \MBF{W_2}&\MBF{I_2}&\MBF{V_2}&& \\
 &\ddots&\ddots&\ddots&\ddots \\
 &&&\MBF{W_p}&\MBF{I_p}
\end{array}\right],
\end{equation}
where $\MBF{I_i}$ denotes an identity matrix of size $n_i$ and $\MBF{D_i}\equiv\MBF{A_i}$.
The $\MBF{V_i}$ and $\MBF{W_i}$ matrices give the SPIKE algorithm its name, because their non-zero elements form tall, narrow submatrices of size $n_i \times k$
(a.k.a. spikes). The equations for these matrices are: 

\begin{equation}
  \label{eq:spikes}
\MBF{V_i} =  \MBF{A_i}^{-1}\MBF{B_i}
;\quad
\MBF{W_i} = \MBF{A_i}^{-1} \MBF{C_i}.
\end{equation}

One source of SPIKE variants is the treatment of the $\MBF{V}$ and $\MBF{W}$ matrices.
In the recursive version of SPIKE that is outlined in this paper, only the bottom $k\times k$ tips of $\MBF{V}$ and $\MBF{W}$ need to be explicitly computed.
Whenever necessary, the forms $\MBF{A_i}^{-1}\MBF{B_i}$ and $\MBF{A_i}^{-1} \MBF{C_i}$ will be used in the place of the corresponding $\MBF{V_i}$ and $\MBF{W_i}$ spikes.

Using the $\MBF{DS}$ on the original problem $\MBF{A}\MBF{X} = \MBF{DS}\MBF{X} = \MBF{F}$, it
can now be broken up into two subproblems, the D stage and the S stage i.e.

\begin{equation}
\MBF{D} \MBF{Y} = 
\left[\begin{array}{ccccc}
\MBF{D_1}&&&& \\
 &\MBF{D_2}&&& \\
 &&\ddots&& \\
 &&&&\MBF{D_p}
\end{array}\right]
\left[
\begin{array}{c}
\MBF{Y_1} \\
\MBF{Y_2} \\
\vdots \\
\MBF{Y_p} \\
\end{array}
\right]
=
\left[
\begin{array}{c}
\MBF{F_1} \\
\MBF{F_2} \\
\vdots \\
\MBF{F_p} \\
\end{array}
\right],
\label{eq:Dstage}
\end{equation}

\begin{equation}
\MBF{S}\MBF{X} =
\left[\begin{array}{ccccc}
\MBF{I_1}&\MBF{V_1}&&& \\
 \MBF{W_2}&\MBF{I_2}&\MBF{V_2}&& \\
 &\ddots&\ddots&\ddots&\ddots \\
 &&&\MBF{W_p}&\MBF{I_p}
\end{array}\right]
\left[
\begin{array}{c}
\MBF{X_1} \\
\MBF{X_2} \\
\vdots \\
\MBF{X_p} \\
\end{array}
\right]
=
\left[
\begin{array}{c}
\MBF{Y_1} \\
\MBF{Y_2} \\
\vdots \\
\MBF{Y_p} \\
\end{array}
\right].\label{eq:Sstage}
\end{equation}

The submatrices of $\MBF{D}$ are decoupled, so the D-stage is straightforward.  
Each partition in (\ref{eq:Dstage}) is solved independently since 
\begin{equation}
\MBF{Y_i} = \MBF{D_i}^{-1} \MBF{F_i}.
\end{equation}

In turn, the vectors and matrices involved in the $\MBF{S}$ stage can be partitioned as follows:

\begin{equation}
\MBF{V_i} = 
\left[\begin{BMAT}(b,18pt,18pt){c}{ccc}
 \MBF{{V_{it}}}  \\ %&\ep{\MBF{0\dots 0}} \\
 \MBF{\tilde{V_i}} \\ %&\ep{\MBF{0\dots 0}}\\
 \MBF{{V_{ib}}}  %&\ep{\MBF{0\dots 0}} 
\end{BMAT}\right]
;\quad
\MBF{W_i} = 
\left[\begin{BMAT}(b,18pt,18pt){c}{ccc}
    %\ep{\MBF{0\dots 0}}&
    \MBF{{W_{it}}}   \\
    %\ep{\MBF{0\dots 0}}&
    \MBF{\tilde{W}_i} \\
    %\ep{\MBF{0\dots 0}}&
    \MBF{{W}_{ib}}   
\end{BMAT}\right],
\end{equation}

\begin{equation}
\MBF{X_i} = 
\left[\begin{BMAT}(b,18pt,18pt){c}{ccc}
 \MBF{{X}_{it}}   \\
 \MBF{\tilde{X}_i} \\
 \MBF{{X}_{ib}}   
\end{BMAT}\right]
;
\quad
\MBF{Y_i} = 
\left[\begin{BMAT}(b,18pt,18pt){c}{ccc}
 \MBF{{Y}_{it}}   \\
 \MBF{\tilde{Y}_i} \\
 \MBF{{Y}_{ib}}   
\end{BMAT}\right],
\end{equation}
where each submatrix denoted with a subscript $t$ or $b$ has a height of $\BW$ rows.
The non-zero partitions of $\MBF{W}_i$ and $\MBF{V}_i$ are $\BW$ columns wide.
Essentially, we have broken out the values coupling the domains of $\MBF{A}$.
Equation  (\ref{eq:Sstage}) can be rewritten as:

\begin{equation}
   \label{eq:solve1}
\left[\begin{BMAT}(b,18pt,18pt){c}{ccc}
\MBF{Y^{}_{1t}} \\
\MBF{\tilde{Y}_1} \\
\MBF{Y^{}_{1b}}
\end{BMAT}\right]
= 
\left[\begin{BMAT}(b,18pt,18pt){c}{ccc}
\MBF{X^{}_{1t}} \\
\MBF{\tilde{X}_1} \\
\MBF{X^{}_{1b}}
\end{BMAT}\right]
+ 
\left[\begin{BMAT}(b,18pt,18pt){c}{ccc}
\MBF{V_{1t}} \\
\MBF{\tilde{V}_1} \\
\MBF{V_{1b}}
\end{BMAT}\right]
\MBF{X^{}_{2t}},
\end{equation}
\begin{equation}
   \label{eq:solve2}
\left[\begin{BMAT}(b,18pt,18pt){c}{ccc}
\MBF{Y_{it}} \\
\MBF{\tilde{Y_i}} \\
\MBF{Y_{ib}}
\end{BMAT}\right]
= 
\left[\begin{BMAT}(b,18pt,18pt){c}{ccc}
\MBF{X_{it}} \\
\MBF{\tilde{X_i}} \\
\MBF{X_{ib}}
\end{BMAT}\right]
+ 
\left[\begin{BMAT}(b,18pt,18pt){c}{ccc}
\MBF{V_{it}} \\
\MBF{\tilde{V_i}} \\
\MBF{V_{ib}}
\end{BMAT}\right]
\MBF{X_{i+1t}} 
+ 
\left[\begin{BMAT}(b,18pt,18pt){c}{ccc}
\MBF{W_{it}} \\
\MBF{\tilde{W_i}} \\
\MBF{W_{ib}}
\end{BMAT}\right]
\MBF{X_{i-1b}},\quad \mbox{for} \quad i \in 2\dots p-1.
\end{equation}
\begin{equation}
   \label{eq:solve3}
\left[\begin{BMAT}(b,18pt,18pt){c}{ccc}
\MBF{Y_{pt}} \\
\MBF{\tilde{Y}_p} \\
\MBF{Y_{pb}}
\end{BMAT}\right]
= 
\left[\begin{BMAT}(b,18pt,18pt){c}{ccc}
\MBF{X_{pt}} \\
\MBF{\tilde{X}_p} \\
\MBF{X_{pb}}
\end{BMAT}\right]
+ 
\left[\begin{BMAT}(b,18pt,18pt){c}{ccc}
\MBF{W_{pt}} \\
\MBF{\tilde{W}_p} \\
\MBF{W_{pb}}
\end{BMAT}\right]
\MBF{X_{p-1b}}.
\end{equation}
Interestingly, the large middle sections of these vectors may be ignored at first.
This will lead to the following definition of the tops and bottoms of these vectors that is amenable to reduced system formation:  

\begin{equation}
\left[\begin{BMAT}(b,18pt,18pt){c}{cc}
\MBF{Y^{}_{1t}} \\
\MBF{Y^{}_{1b}}
\end{BMAT}\right]
= 
\left[\begin{BMAT}(b,18pt,18pt){c}{cc}
\MBF{X^{}_{1t}} \\
\MBF{X^{}_{1b}}
\end{BMAT}\right]
+ 
\left[\begin{BMAT}(b,18pt,18pt){c}{cc}
\MBF{V_{1t}} \\
\MBF{V_{1b}}
\end{BMAT}\right]
\MBF{X^{}_{2t}}, 
\end{equation}
\begin{equation}
\left[\begin{BMAT}(b,18pt,18pt){c}{cc}
\MBF{Y_{it}} \\
\MBF{Y_{ib}}
\end{BMAT}\right]
= 
\left[\begin{BMAT}(b,18pt,18pt){c}{cc}
\MBF{X_{it}} \\
\MBF{X_{ib}}
\end{BMAT}\right]
+ 
\left[\begin{BMAT}(b,18pt,18pt){c}{cc}
\MBF{V_{it}} \\
\MBF{V_{ib}}
\end{BMAT}\right]
\MBF{X_{i+1t}} 
+ 
\left[\begin{BMAT}(b,18pt,18pt){c}{cc}
\MBF{W_{it}} \\
\MBF{W_{ib}}
\end{BMAT}\right]
\MBF{X_{i-1b}},\quad \mbox{for} \quad i \in 2\dots p-1.
\end{equation}
\begin{equation}
\left[\begin{BMAT}(b,18pt,18pt){c}{cc}
\MBF{Y_{pt}} \\
\MBF{Y_{pb}}
\end{BMAT}\right]
= 
\left[\begin{BMAT}(b,18pt,18pt){c}{cc}
\MBF{X_{pt} }\\
\MBF{X_{pb}}
\end{BMAT}\right]
+ 
\left[\begin{BMAT}(b,18pt,18pt){c}{cc}
\MBF{W_{pt}} \\
\MBF{W_{pb}}
\end{BMAT}\right]
\MBF{X_{p-1b}}.
\end{equation}

\begin{figure}[htbp] 
\includegraphics[keepaspectratio,width=\textwidth]{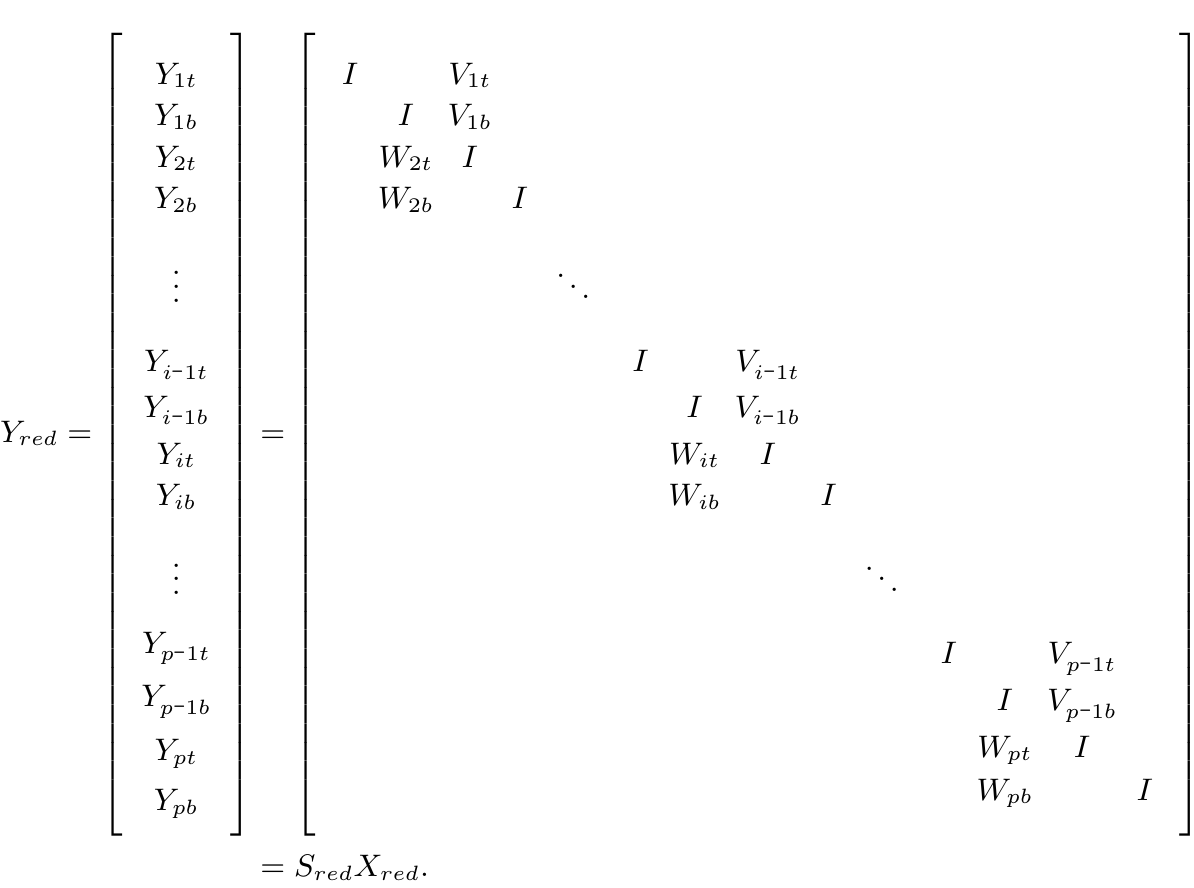} 
\caption{The SPIKE reduced system.}
\label{fig:red}
\end{figure}

The reduced system is shown in Figure~\ref{fig:red}. Conceptually, the reduced system could just be thought of as a small banded matrix problem.
One common source of SPIKE variants is the specific method of solving this reduced system.
The `recursive method' for solving the reduced system is discussed in the next section.

Once the reduced system is solved, we obtain the values for $\MBF{X_{b,i}}$ and $\MBF{X_{t,i}}$ with $i \in 1\dots p$.
In turn, the values for $\MBF{\tilde{X}_{i}}$ for all $i$ can be straightforwardly recovered using (\ref{eq:solve1}), (\ref{eq:solve2}), and (\ref{eq:solve3})
  (a.k.a. the retrieval stage). In some practical implementations of SPIKE,
  once the factorization stage is complete, the middle part of the spikes $V$ and $W$ (resp. $\tilde{V}$ and $\tilde{W}$)
  are not stored in memory, so they are not available during the retrieval stage.
  In addition, we note that the spikes $V_1$ and $W_p$ are never explicitly computed providing further optimization of the algorithm (cf. section~\ref{opt_sweeps}).
  Consequently, the spikes can instead be replaced by their expression in (\ref{eq:spikes}) leading to
  the following solve operations:

\begin{equation} \label{eq:solve1n}
  X_1=Y_1-A_1^{-1}({B}_1\MBF{X_{2t}}),
\end{equation}

\begin{equation} \label{eq:solve2n}
  X_i=Y_i-A_i^{-1}({B}_i\MBF{X_{i+1t}}+{C}_i\MBF{X_{i-1b}}),\quad \mbox{for} \quad i \in 2\dots p-1,
\end{equation}

\begin{equation} \label{eq:solve3n}
  X_p=Y_p-A_p^{-1}({C}_p\MBF{X_{p-1b}}).
\end{equation}

At this point, $\MBF{X}$ has been found and the computation is complete.

\subsection{Recursive reduced system} \label{recursive_reduced_section}

The reduced system in Figure~\ref{fig:red}, represents the inter-domain relationships for the partitioning performed on $\MBF{A}$, it is
  of size of $2p\BW$ which scales  
linearly with the number of partitions $p$. In order to fully capitalize on the 
  performance gained by exploiting parallelism in the factorization and retrieval stages, the reduced system should not be explicitly formed.
Among the multiple techniques that are available for solving the reduced system in parallel, the recursive SPIKE technique provides
the best trade-off between generality and parallel efficiency.
A full derivation of the recursive method for solving the reduced system is shown in \cite{Polizzi:2006}.
The essential observation is that the reduced system is banded, and, as a result, SPIKE may be used to solve it. From the original reduced system, a
new spike matrix $S$ will then be generated which, in turn, could be solved by SPIKE with half the number of partitions. The process can be repeated recursively,
where the number of partitions to consider is divided by two at each recursion level,
and  until only two partitions are left.

For clarity, an extra superscript indexing has been added to all the submatrices in the following equations to designate the level of recursion.
Here, the process will be illustrated using a four-partition example (i.e. $p=4$)
which is sufficient to provide one level of recursion and
show the central concept of the scheme. Our starting point is the original
  four-partition reduced system:

\begin{equation}
\label{RED_SYS_RECA}
\MBF{Y^{[1]}}=
\left[\begin{BMAT}(r,12pt,12pt){c}{cccccccccc}
\MBF{Y^{[1]}_{1t}}\\
\MBF{Y^{[1]}_{1b}}\\
\MBF{Y^{[1]}_{2t}}\\
\MBF{Y^{[1]}_{2b}}\\
\MBF{Y^{[1]}_{3t}}\\
\MBF{Y^{[1]}_{3b}}\\
\MBF{Y^{[1]}_{4t}}\\
\MBF{Y^{[1]}_{4b}}\\
\end{BMAT} \right] 
=
\left[\begin{BMAT}(m,12pt,12pt){cccccccc}{cccccccc}
 \MBF{I}&              &\MBF{V^{[1]}_{1t}}&              &              &              &              &       \\
        &\MBF{I}       &\MBF{V^{[1]}_{1b}}&              &              &              &              &       \\
        &\MBF{W^{[1]}_{2t}}&\MBF{I}       &              &\MBF{V^{[1]}_{2t}}&              &              &       \\
        &\MBF{W^{[1]}_{2b}}&              &\MBF{I}       &\MBF{V^{[1]}_{2b}}&              &              &       \\
        &              &              &\MBF{W^{[1]}_{3t}}&\MBF{I}       &              &\MBF{V^{[1]}_{3t}}&       \\
        &              &              &\MBF{W^{[1]}_{3b}}&              &\MBF{I}       &\MBF{V^{[1]}_{3b}}&       \\
        &              &              &              &              &\MBF{W^{[1]}_{4t}}&\MBF{I}       &       \\
        &              &              &              &              &\MBF{W^{[1]}_{4b}}&              &\MBF{I}\\
\end{BMAT}\right]
\left[\begin{BMAT}(r,12pt,12pt){c}{cccccccccc}
\MBF{X^{[1]}_{1t}}\\
\MBF{X^{[1]}_{1b}}\\
\MBF{X^{[1]}_{2t}}\\
\MBF{X^{[1]}_{2b}}\\
\MBF{X^{[1]}_{3t}}\\
\MBF{X^{[1]}_{3b}}\\
\MBF{X^{[1]}_{4t}}\\
\MBF{X^{[1]}_{4b}}\\
\end{BMAT} \right]
=
S^{[1]}\MBF{X^{[1]}},
\end{equation}

where we use the notation $Y^{[1]}=Y_{red}$, $S^{[1]}=S_{red}$ and $X^{[1]}=X_{red}$ to
  emphasize the current level of recursion (level one here).
We then perform a new $DS$ SPIKE factorization of the reduced system using half the number of  partitions (so two partitions here),

\begin{equation}
\label{RED_SYS_RECB}
  \MBF{S^{[1]}}
  =
\left[\begin{BMAT}(m,10pt,12pt){cccc:cccc}{cccc:cccc}
 \MBF{I}&              &\MBF{V^{[1]}_{1t}}&              &              &              &              &       \\
        &\MBF{I}       &\MBF{V^{[1]}_{1b}}&              &              &              &              &       \\
        &\MBF{W^{[1]}_{2t}}&\MBF{I}       &              &              &              &              &       \\
        &\MBF{W^{[1]}_{2b}}&              &\MBF{I}       &              &              &              &       \\
        &              &              &\             &\MBF{I}       &              &\MBF{V^{[1]}_{3t}}&       \\
        &              &              &\             &              &\MBF{I}       &\MBF{V^{[1]}_{3b}}&       \\
        &              &              &              &              &\MBF{W^{[1]}_{4t}}&\MBF{I}       &       \\
        &              &              &              &              &\MBF{W^{[1]}_{4b}}&              &\MBF{I}\\
\end{BMAT}\right]
\left[\begin{BMAT}(m,10pt,12pt){cccc:cccc}{cccc:cccc}
 \MBF{I}&              &              &              &\MBF{V^{[2]}_{1t}}&              &              &       \\
        &\MBF{I}       &              &              &\MBF{V^{[2]}_{1b}}&              &              &       \\
        &              &\MBF{I}       &              &\MBF{V^{[2]}_{2t}}&              &              &       \\
        &              &              &\MBF{I}       &\MBF{V^{[2]}_{2b}}&              &              &       \\
        &              &              &\MBF{W^{[2]}_{3t}}&\MBF{I}       &              &              &       \\
        &              &              &\MBF{W^{[2]}_{3b}}&              &\MBF{I}       &              &       \\
        &              &              &\MBF{W^{[2]}_{4t}}&              &              &\MBF{I}       &       \\
        &              &              &\MBF{W^{[2]}_{4b}}&              &              &              &\MBF{I}\\
\end{BMAT}\right] =\MBF{D^{[1]}S^{[2]}},
\end{equation}

with
\begin{equation}
\left[\begin{BMAT}(r,12pt,12pt){cccc}{cccc}
\MBF{I}&              &\MBF{V^{[1]}_{1t}}&         \\     
       &\MBF{I}       &\MBF{V^{[1]}_{1b}}&         \\     
       &\MBF{W^{[1]}_{2t}}&\MBF{I}       &         \\     
       &\MBF{W^{[1]}_{2b}}&              &\MBF{I}  \\     
\end{BMAT}\right]
\left[\begin{BMAT}(r,12pt,12pt){c}{cccc}
 \MBF{V^{[2]}_{1t}}\\
 \MBF{V^{[2]}_{1b}}\\
 \MBF{V^{[2]}_{2t}}\\
 \MBF{V^{[2]}_{2b}}\\
\end{BMAT}\right]
=
\left[\begin{BMAT}(r,12pt,12pt){c}{cccc}
 \MBF{   0    }\\
 \MBF{   0    }\\
 \MBF{V^{[1]}_{2t}}\\
 \MBF{V^{[1]}_{2b}}\\
\end{BMAT}\right]
\rightarrow
\begin{cases}
&\left[\begin{BMAT}(r,12pt,12pt){cc}{cc}
\MBF{I}       &\MBF{V^{[1]}_{1b}}\\     
\MBF{W^{[1]}_{2t}}&\MBF{I}       \\     
\end{BMAT}\right]
\left[\begin{BMAT}(r,12pt,12pt){c}{cccc}
 \MBF{V^{[2]}_{1b}}\\
 \MBF{V^{[2]}_{2t}}\\
\end{BMAT}\right]
=
\left[\begin{BMAT}(r,12pt,12pt){c}{cccc}
 \MBF{    0   }\\
 \MBF{V^{[1]}_{2t}}\\
\end{BMAT}\right]\\[15pt]
&\MBF{V^{[2]}_{1t} = - V^{[1]}_{1t} V^{[2]}_{2t}}\\[5pt]
&\MBF{V^{[2]}_{2b} = V^{[1]}_{2b} - W^{[1]}_{2b} V^{[2]}_{1b}}\\
\end{cases}
\end{equation} 

and

\begin{equation}
\label{RED_SYS_W}
\left[\begin{BMAT}(r,12pt,12pt){cccc}{cccc}
\MBF{I}&              &\MBF{V^{[1]}_{3t}}&         \\     
       &\MBF{I}       &\MBF{V^{[1]}_{3b}}&         \\  
       &\MBF{W^{[1]}_{4t}}&\MBF{I}       &         \\     
       &\MBF{W^{[1]}_{4b}}&              &\MBF{I}  \\     
\end{BMAT}\right]
\left[\begin{BMAT}(r,12pt,12pt){c}{cccc}
 \MBF{W^{[2]}_{3t}}\\
 \MBF{W^{[2]}_{3b}}\\
 \MBF{W^{[2]}_{4t}}\\
 \MBF{W^{[2]}_{4b}}\\
\end{BMAT}\right]
=
\left[\begin{BMAT}(r,12pt,12pt){c}{cccc}
 \MBF{W^{[1]}_{3t}}\\
 \MBF{W^{[1]}_{3t}}\\
 \MBF{  0     }\\
 \MBF{  0     }\\
\end{BMAT}\right]
\rightarrow
\begin{cases}
&
\left[\begin{BMAT}(r,12pt,12pt){cc}{cc}
\MBF{I}       &\MBF{V^{[1]}_{3b}}\\     
\MBF{W^{[1]}_{4t}}&\MBF{I}       \\     
\end{BMAT}\right]
\left[\begin{BMAT}(r,12pt,12pt){c}{cccc}
 \MBF{W^{[2]}_{3b}}\\
 \MBF{W^{[2]}_{4t}}\\
\end{BMAT}\right]
=
\left[\begin{BMAT}(r,12pt,12pt){c}{cccc}
 \MBF{W^{[1]}_{3t}}\\
 \MBF{  0     }\\
\end{BMAT}\right]\\[15pt]
&\MBF{W^{[2]}_{3t} = W^{[1]}_{3t} - V^{[1]}_{3t} W^{[2]}_{4t}}\\[5pt]
&\MBF{W^{[2]}_{4b} = - W^{[1]}_{4b} W^{[2]}_{3b}}\\
\end{cases}
\end{equation}

It should be noted that the widths of the $V$ and $W$ spikes in $\MBF{S^{[2]}}$ are equal to the widths of $\MBF{V^{[1]}_{2}}$ and $\MBF{W^{[1]}_{3}}$ respectively.
The matrix $\MBF{S^{[2]}}$ is 
already in the form of a two-partition S-matrix,
so the recursion stops at this step.  The reduced system factorization is then complete. 
Solving the reduced system (\ref{RED_SYS_RECB}) can be performed in two stages: (i) Obtain the intermediate solution $Y^{[2]}$

\begin{equation} \label{D1Inverse}
\MBF{D^{[1]} Y^{[2]} = Y^{[1]}}, 
\end{equation}
and (ii) Solve for $X^{[1]}$
\begin{equation} \label{S2Inverse}
\MBF{S^{[2]} X^{[1]} = Y^{[2]}}. 
\end{equation}

First we will look at equation~(\ref{D1Inverse}). The blocks of the $\MBF{D^{[1]}}$ matrix are uncoupled, so they can be solved in parallel. 
In addition, the individual blocks take a form similar to that of a two-partition S-matrix, so an even smaller reduced system can be extracted from each. 

%\small
\begin{equation} \label{D1InverseBlock}
\left[\begin{BMAT}(m,16pt,14pt){cccc:cccc}{cccc:cccc}
 \MBF{I}&              &\MBF{V^{[1]}_{1t}}&              &              &              &              &       \\
        &\MBF{I}       &\MBF{V^{[1]}_{1b}}&              &              &              &              &       \\
        &\MBF{W^{[1]}_{2t}}&\MBF{I}       &              &              &              &              &       \\
        &\MBF{W^{[1]}_{2b}}&              &\MBF{I}       &              &              &              &       \\
        &              &              &\             &\MBF{I}       &              &\MBF{V^{[1]}_{3t}}&       \\
        &              &              &\             &              &\MBF{I}       &\MBF{V^{[1]}_{3b}}&       \\
        &              &              &              &              &\MBF{W^{[1]}_{4t}}&\MBF{I}       &       \\
        &              &              &              &              &\MBF{W^{[1]}_{4b}}&              &\MBF{I}\\
\end{BMAT}\right]
\left[\begin{BMAT}(r,12pt,14pt){c}{cccccccc}
\MBF{Y^{[2]}_{1t}}\\
\MBF{Y^{[2]}_{1b}}\\
\MBF{Y^{[2]}_{2t}}\\
\MBF{Y^{[2]}_{2b}}\\
\MBF{Y^{[2]}_{3t}}\\
\MBF{Y^{[2]}_{3b}}\\
\MBF{Y^{[2]}_{4t}}\\
\MBF{Y^{[2]}_{4b}}\\
\end{BMAT} \right]
=
\left[\begin{BMAT}(r,12pt,14pt){c}{cccccccc}
\MBF{Y^{[1]}_{1t}}\\
\MBF{Y^{[1]}_{1b}}\\
\MBF{Y^{[1]}_{2t}}\\
\MBF{Y^{[1]}_{2b}}\\
\MBF{Y^{[1]}_{3t}}\\
\MBF{Y^{[1]}_{3b}}\\
\MBF{Y^{[1]}_{4t}}\\
\MBF{Y^{[1]}_{4b}}\\
\end{BMAT} \right]
\end{equation}
%\normalsize

%\small
\begin{equation}
\left[\begin{BMAT}(r,12pt,12pt){cccc}{cccc}
\MBF{I}&              &\MBF{V^{[1]}_{1t}}&         \\     
       &\MBF{I}       &\MBF{V^{[1]}_{1b}}&         \\     
       &\MBF{W^{[1]}_{2t}}&\MBF{I}       &         \\     
       &\MBF{W^{[1]}_{2b}}&              &\MBF{I}  \\     
\end{BMAT}\right]
\left[\begin{BMAT}(r,12pt,12pt){c}{cccc}
 \MBF{Y^{[2]}_{1t}}\\
 \MBF{Y^{[2]}_{1b}}\\
 \MBF{Y^{[2]}_{2t}}\\
 \MBF{Y^{[2]}_{2b}}\\
\end{BMAT}\right]
=
\left[\begin{BMAT}(r,12pt,12pt){c}{cccc}
 \MBF{Y^{[1]}_{1t}}\\
 \MBF{Y^{[1]}_{1t}}\\
 \MBF{Y^{[1]}_{2t}}\\
 \MBF{Y^{[1]}_{2b}}\\
\end{BMAT}\right]
\rightarrow
\begin{cases}
&\left[\begin{BMAT}(r,12pt,12pt){cc}{cc}
\MBF{I}       &\MBF{V^{[1]}_{1b}}\\     
\MBF{W^{[1]}_{2t}}&\MBF{I}       \\     
\end{BMAT}\right]
\left[\begin{BMAT}(r,12pt,12pt){c}{cccc}
 \MBF{Y^{[2]}_{1b}}\\
 \MBF{Y^{[2]}_{2t}}\\
\end{BMAT}\right]
=
\left[\begin{BMAT}(r,12pt,12pt){c}{cccc}
 \MBF{Y^{[1]}_{1t}}\\
 \MBF{Y^{[1]}_{2t}}\\
\end{BMAT}\right]\\
&\MBF{Y^{[2]}_{1t} = Y^{[1]}_{1t} - V^{[1]}_{1t} Y^{[2]}_{2t}}\\
&\MBF{Y^{[2]}_{2b} = Y^{[1]}_{2b} - W^{[1]}_{2b} Y^{[2]}_{1b}}\\
\end{cases}
\end{equation} 

\begin{equation}
\left[\begin{BMAT}(r,12pt,12pt){cccc}{cccc}
\MBF{I}&              &\MBF{V^{[1]}_{3t}}&         \\     
       &\MBF{I}       &\MBF{V^{[1]}_{3b}}&         \\  
       &\MBF{W^{[1]}_{4t}}&\MBF{I}       &         \\     
       &\MBF{W^{[1]}_{4b}}&              &\MBF{I}  \\     
\end{BMAT}\right]
\left[\begin{BMAT}(r,12pt,12pt){c}{cccc}
 \MBF{Y^{[2]}_{3t}}\\
 \MBF{Y^{[2]}_{3b}}\\
 \MBF{Y^{[2]}_{4t}}\\
 \MBF{Y^{[2]}_{4b}}\\
\end{BMAT}\right]
=
\left[\begin{BMAT}(r,12pt,12pt){c}{cccc}
 \MBF{Y^{[1]}_{3t}}\\
 \MBF{Y^{[1]}_{3t}}\\
 \MBF{Y^{[1]}_{4t}}\\
 \MBF{Y^{[1]}_{4b}}\\
\end{BMAT}\right]
\rightarrow
\begin{cases}
&
\left[\begin{BMAT}(r,12pt,12pt){cc}{cc}
\MBF{I}       &\MBF{V^{[1]}_{3b}}\\     
\MBF{W^{[1]}_{4t}}&\MBF{I}       \\     
\end{BMAT}\right]
\left[\begin{BMAT}(r,12pt,12pt){c}{cccc}
 \MBF{Y^{[2]}_{3b}}\\
 \MBF{Y^{[2]}_{4t}}\\
\end{BMAT}\right]
=
\left[\begin{BMAT}(r,12pt,12pt){c}{cccc}
 \MBF{Y^{[1]}_{3t}}\\
 \MBF{Y^{[1]}_{4t}}\\
\end{BMAT}\right]\\
&\MBF{Y^{[2]}_{3t} = Y^{[1]}_{3t} - V^{[1]}_{3t} Y^{[2]}_{4t}}\\
&\MBF{Y^{[2]}_{4b} = Y^{[1]}_{4b} - W^{[1]}_{4b} Y^{[2]}_{3b}}\\
\end{cases}
\end{equation}

Therefore, the ${\MBF{D_1}}$ matrix solve has been reduced to two $2\BW \times 2\BW$ solve operations, which are performed in parallel, and some recovery operations. 
Next, equation~(\ref{S2Inverse}) must be solved. 
This is simply a two-partition S-matrix, so we will extract a reduced system and perform recovery sweeps as usual,

\begin{equation}
\left[\begin{BMAT}(m,16pt,14pt){cccc:cccc}{cccc:cccc}
 \MBF{I}&              &              &              &\MBF{V^{[2]}_{1t}}&              &              &       \\
        &\MBF{I}       &              &              &\MBF{V^{[2]}_{1b}}&              &              &       \\
        &              &\MBF{I}       &              &\MBF{V^{[2]}_{2t}}&              &              &       \\
        &              &              &\MBF{I}       &\MBF{V^{[2]}_{2b}}&              &              &       \\
        &              &              &\MBF{W^{[2]}_{3t}}&\MBF{I}       &              &              &       \\
        &              &              &\MBF{W^{[2]}_{3b}}&              &\MBF{I}       &              &       \\
        &              &              &\MBF{W^{[2]}_{4t}}&              &              &\MBF{I}       &       \\
        &              &              &\MBF{W^{[2]}_{4b}}&              &              &              &\MBF{I}\\
\end{BMAT}\right]
\left[\begin{BMAT}(r,14pt,14pt){c}{cccccccc}
\MBF{X^{[1]}_{1t}}\\
\MBF{X^{[1]}_{1b}}\\
\MBF{X^{[1]}_{2t}}\\
\MBF{X^{[1]}_{2b}}\\
\MBF{X^{[1]}_{3t}}\\
\MBF{X^{[1]}_{3b}}\\
\MBF{X^{[1]}_{4t}}\\
\MBF{X^{[1]}_{4b}}\\
\end{BMAT} \right]
=
\left[\begin{BMAT}(r,14pt,14pt){c}{cccccccc}
\MBF{Y^{[2]}_{1t}}\\
\MBF{Y^{[2]}_{1b}}\\
\MBF{Y^{[2]}_{2t}}\\
\MBF{Y^{[2]}_{2b}}\\
\MBF{Y^{[2]}_{3t}}\\
\MBF{Y^{[2]}_{3b}}\\
\MBF{Y^{[2]}_{4t}}\\
\MBF{Y^{[2]}_{4b}}\\
\end{BMAT} \right]
\end{equation}

\begin{equation}
\begin{split}
\left[\begin{BMAT}(r,12pt,12pt){cc}{cc}
\MBF{I}       &\MBF{V^{[2]}_{2b}}\\     
\MBF{W^{[2]}_{3t}}&\MBF{I}       \\     
\end{BMAT}\right]
\left[\begin{BMAT}(r,12pt,12pt){c}{cccc}
 \MBF{X^{[1]}_{2b}}\\
 \MBF{X^{[1]}_{3t}}\\
\end{BMAT}\right]
=
\left[\begin{BMAT}(r,12pt,12pt){c}{cccc}
 \MBF{Y^{[2]}_{2t}}\\
 \MBF{Y^{[2]}_{3t}}\\
\end{BMAT}\right]\\
\left[
\begin{BMAT}(r,12pt,12pt){c}{ccc}
 \MBF{X^{[1]}_{1t}} \\
 \MBF{X^{[1]}_{1b}} \\
 \MBF{X^{[1]}_{2t}} 
\end{BMAT}
\right]
=
\left[
\begin{BMAT}(r,12pt,12pt){c}{ccc}
 \MBF{Y^{[2]}_{1t}} \\
 \MBF{Y^{[2]}_{1b}} \\
 \MBF{Y^{[2]}_{2t}} 
\end{BMAT}
\right]
-
\left[
\begin{BMAT}(r,12pt,12pt){c}{ccc}
 \MBF{V^{[2]}_{1t}} \\
 \MBF{V^{[2]}_{1b}} \\
 \MBF{V^{[2]}_{2t}} 
\end{BMAT}
\right]
\MBF{X^{[1]}_{3t}} \\
\left[
\begin{BMAT}(r,12pt,12pt){c}{ccc}
 \MBF{X^{[1]}_{3b}} \\
 \MBF{X^{[1]}_{4t}} \\
 \MBF{X^{[1]}_{4b}} 
\end{BMAT}
\right]
=
\left[
\begin{BMAT}(r,12pt,12pt){c}{ccc}
 \MBF{Y^{[2]}_{3b}} \\
 \MBF{Y^{[2]}_{4t}} \\
 \MBF{Y^{[2]}_{4b}} 
\end{BMAT}
\right]
-
\left[
\begin{BMAT}(r,12pt,12pt){c}{ccc}
 \MBF{W^{[2]}_{3b}} \\
 \MBF{W^{[2]}_{4t}} \\
 \MBF{W^{[2]}_{4b}} 
\end{BMAT}
\right]
\MBF{X^{[1]}_{2b}}
\end{split}
\end{equation}

At this point the $\MBF{X^{[1]}}$ vectors have been found, so the reduced system is solved.
The total number of $2\BW \times 2k$ solve operations is the same as the number of partition interfaces, $p-1$.
The total computational cost spent on solve operations is $O(p \times \BW \times n_{rhs})$.
However, all the solve operations in each recursive level may be performed in parallel.
Because the system is split in half with each recursive level, the total number of recursive levels is $\log_2(p)$.
Therefore, the combined critical path length of all the solve operations in the solve stage is $O(\log_2(p)\times \BW \times n_{rhs})$.
For the same reason, the reduced system factorization stage solve operations have a critical path length of $O(\log_2(p)\times \BW^2 )$.
So, the total cost of the solve operations is $O(\log_2(p)\times \BW \times \max(k,n_{rhs}))$.
There is also some overhead involved with the solution recovery operations and communication, but this has not been found to be significant. 

This completes the description of the recursive reduced system.
This method of solving the reduced system can significantly improve performance by exploiting parallelism in the problem. 
However, because the procedure progresses through recursive levels by repeatedly splitting submatrices in half, the recursive reduced system
limits the number of partitions allowable to a power of two. 
A method of decoupling the number of threads used from the number of partitions will be shown in Section \ref{Increased_parallelism_for_Recursive_Spike}.
Next, we look at optimizations specific to the banded structure. 

\subsection{Optimizing per-partition costs} \label{opt_sweeps}

In Section~\ref{basics_section}, we neglected the specifics of
the factorization performed on the blocks, $\MBF{D}_i$.
The primary computational costs for SPIKE are the matrix operations performed on each block. 
The goal, then, is to reduce the number of solve operations performed.

The $\MBF{D_i}$ matrices are factorized into triangular matrices. 
For a total number of partitions $p$, partitions $1$ to $p-1$ use an LU factorization. 
For the final partition, a UL factorization is used. 
In practice, non-pivoting factorizations have been used to retain the pattern of zeroes in the $\MBF{B_i}$ and $\MBF{C_i}$ matrices. 
In Section~\ref{Pivoting_Spike_Section} a method of overcoming this limitation and applying partial pivoting will be shown. 
In the following, we will be working with the non-pivoting SPIKE algorithm
  using the diagonal boosting strategy originally introduced in \cite{Polizzi:2006}
  that offers an excellent
trade-off between accuracy and performance.
The first detail to look at is the creation of the V spikes,

\begin{equation} \label{V_formula_end}
\MBF{V_i} = \MBF{A_i^{-1}} \MBF{B_i}
=
\MBF{U_i^{-1}}\MBF{L_i^{-1}}
\left[ \begin{BMAT}(b,18pt,18pt){c}{bc:c}
\matmid{\MBF{0}}\\ 
 \MBF{\hat{B}_i}
  \end{BMAT} \right].
%\left[\begin{BMAT}(b,18pt,18pt){c}{cc}
% \MBF{0} \\ % &\MBF{0} \\
% \MBF{\hat{B}_i} % &\MBF{0} 
%\end{BMAT}\right].
\end{equation}

The matrix $\MBF{L}_i^{-1}$ is lower triangular. 
The solve operation for a lower triangular matrix begins by identifying the topmost rows in the solution vector, and works downward. 
For this reason we label this a ``downward sweep''. 
In the case of equation~(\ref{V_formula_end}), the downward sweep is simply passing over zeroes until the topmost rows of $\MBF{\hat{B}_i}$ are reached. 
So, this sweep may be shortened by beginning it at that point. 
This shortens the downward sweep from a height of $n_i$ to a height of $\BW$, rendering it relatively inconsequential in terms of computational cost. 

For the final partition, the matrix is UL factorized.
The optimization is similar, but it instead avoids the zeroes in the upward sweep. 

\begin{equation} \label{W_formula_end}
\MBF{W_p} = \MBF{A_p^{-1}} \MBF{C_p}
=
\MBF{L_p^{-1}}
\MBF{U_p^{-1}}
\left[ \begin{BMAT}(b,18pt,18pt){c}{b:cc}
     \MBF{\hat{C}_p} \\
\matmid{\MBF{0}}
  \end{BMAT} \right],
%\left[\begin{BMAT}(b,18pt,18pt){c}{cc}
%    %\MBF{0}  &
%    \MBF{\hat{C}_p} \\
%    %\MBF{0}  &
%    \MBF{0} 
%\end{BMAT}\right].
\end{equation}

The next important variation from the basic version of SPIKE discussed earlier is the treatment of the V and W spikes. 
Using the definitions for $\MBF{V}_i$ and $\MBF{W}_i$ above, and the fact that $\MBF{Y}_i=\MBF{D_i^{-1}} \MBF{F_i}$, we may rewrite the retrieval
stage shown previously in (\ref{eq:solve1n}), (\ref{eq:solve2n}), and (\ref{eq:solve3n}), as follows:
%\begin{equation}
 % \MBF{Y_i} = \MBF{D_i^{-1}} \MBF{F_i},
%\end{equation}

\begingroup
\thinmuskip=0mu
\medmuskip=0mu
\thickmuskip=0mu
\small
\begin{equation} \label{one}
\begin{split}
\left[\begin{BMAT}(b,18pt,18pt){c}{ccc}
\MBF{X_{1t}} \\
\MBF{\tilde{X}_1} \\
\MBF{X_{1b}}
\end{BMAT}\right]
=
\left[\begin{BMAT}(b,18pt,18pt){c}{ccc}
\MBF{Y_{1t}} \\
\MBF{\tilde{Y}_1} \\
\MBF{Y_{1b}}
\end{BMAT}\right]
- 
\MBF{V_{1} }
\MBF{X_{2t}}
=
\MBF{A_1^{-1}}
\left(
\left[\begin{BMAT}(b,18pt,18pt){c}{cc:c}
\MBF{F_{1t}} \\
\MBF{\tilde{F}_1} \\
\MBF{F_{1b}}
\end{BMAT}\right]
- 
\left[ 
\begin{BMAT}(b,18pt,18pt){c}{bc:c}
\matmid{\MBF{0}} \\
\MBF{\hat{B}_1}
\end{BMAT} \right]
\MBF{X_{2t}} \right)
=
\MBF{U_1^{-1}}
\left(
\MBF{L_1^{-1}}
\left[\begin{BMAT}(b,18pt,18pt){c}{cc:c}
\MBF{F_{1t}} \\
\MBF{\tilde{F}_1} \\
\MBF{F_{1b}}
\end{BMAT}\right]
- 
\MBF{L_1^{-1}}
\left[ \begin{BMAT}(b,18pt,18pt){c}{bc:c}
\matmid{\MBF{0}}\\ 
 \MBF{\hat{B}_1}
\end{BMAT} \right]
\MBF{X_{2t}} \right)
\end{split},
\end{equation}

\begin{equation}\label{two}
\begin{split}
\left[\begin{BMAT}(b,18pt,18pt){c}{ccc}
\MBF{X_{it}} \\
\MBF{\tilde{X}_i} \\
\MBF{X_{ib}}
\end{BMAT}\right]
= 
\left[\begin{BMAT}(b,18pt,18pt){c}{ccc}
\MBF{Y_{it}} \\
\MBF{\tilde{Y}_i} \\
\MBF{Y_{ib}}
\end{BMAT}\right]
- 
\MBF{V_{i}}
\MBF{X_{i+1t}} 
- 
\MBF{W_{i}} 
\MBF{X_{i-1b}}
=
\MBF{A_i^{-1}}
\left[\begin{BMAT}(b,18pt,18pt){c}{c:c:c}
\MBF{F_{it}} \\
\MBF{\tilde{F}_i} \\
\MBF{F_{ib}}
\end{BMAT}\right]
- 
\MBF{A_i^{-1}}
\left(
\left[\begin{BMAT}(b,18pt,18pt){c}{bc:c}
\matmid{\MBF{0}}\\
\MBF{\hat{B}_i} 
\end{BMAT} \right]
\MBF{X_{i+1t}} 
+ 
\left[\begin{BMAT}(b,18pt,18pt){c}{c:bc}
\MBF{\hat{C}_i} \\
\matmidt{\MBF{0}}
\end{BMAT} \right]
\MBF{X_{i-1b}}
\right)
\end{split},
\end{equation}

\begin{equation}\label{three}
\begin{split}
\left[\begin{BMAT}(b,18pt,18pt){c}{ccc}
\MBF{X_{pt}} \\
\MBF{\tilde{X}_p} \\
\MBF{X_{pb}}
\end{BMAT}\right]
= 
\left[\begin{BMAT}(b,18pt,18pt){c}{ccc}
\MBF{Y_{pt}} \\
\MBF{\tilde{Y}_p} \\
\MBF{Y_{pb}}
\end{BMAT}\right]
- 
\MBF{W_{p}} 
\MBF{X_{p-1b}}
=
\MBF{A_p^{-1}}
\left(
\left[\begin{BMAT}(b,18pt,18pt){c}{c:cc}
\MBF{F_{pt}} \\
\MBF{\tilde{F}_p} \\
\MBF{F_{pb}}
\end{BMAT}\right]
- 
\left[\begin{BMAT}(b,18pt,18pt){c}{c:bc}
\MBF{\hat{C}_p} \\
\matmidt{\MBF{0}}
\end{BMAT} \right]
\MBF{X_{p-1b}}
\right)
=
\MBF{L_p^{-1}}
\left(
\MBF{U_p^{-1}}
\left[\begin{BMAT}(b,18pt,18pt){c}{c:cc}
\MBF{F_{pt}} \\
\MBF{\tilde{F}_p} \\
\MBF{F_{pb}}
\end{BMAT}\right]
- 
\MBF{U_p^{-1}}
\left[\begin{BMAT}(b,18pt,18pt){c}{c:bc}
\MBF{\hat{C}_p} \\
\matmidt{\MBF{0}}
\end{BMAT} \right]
\MBF{X_{p-1b}}
\right)
\end{split}.
\end{equation}

\endgroup

\normalsize

For the first partition, the task of the D stage is to create the bottom tip of the vector $\MBF{A_1}^{-1}\MBF{F}_1$.
Since that vector is unmodified by the reduced system until we reach the very bottom, the L sweep is uncontaminated until it hits that point.
So, we may save the large L sweep from the D stage, and use a small U sweep over the bottom tip to generate the needed values of $\MBF{A_1}^{-1}\MBF{F_1}$.
Similarly, for the last partition, we sweep up across values uncontaminated by the reduced system until we hit the very top of the $\MBF{Y_{p}}$ vector.
In this way, the solve stage for the first and last partitions is performed with just two large sweeps, and a collection of small sweeps and multiplications with practically no cost. 
For all other partitions, a total of four sweeps per partition are needed in the solve stage. 
	
The reduced system only needs $\MBF{V_{1b}}$ for the first partition, and $\MBF{W_{pt}}$ for the last partition.
As a result the upward sweep in equation~(\ref{V_formula_end}) can also be truncated. 
Similarly, the downward sweep in equation~(\ref{W_formula_end}) is truncated.
This results in no full sweeps in these partitions during the factorization stage.
For the middle partitions, the tips of $V$ and $W$ can be obtained using three full sweeps in the SPIKE factorization stage,
one full sweep to generate the spike $V$ and two full sweeps to generate $W$.

\begin{table}[htbp]
  \begin{center}
    \begin{tabular}{|l||c|c|}
      \hline
  \# of full sweeps & Factorization stage & Solve stage \\ \hline
  First \& Last partition & 0 & 2 \\  \hline
  Middle partitions &3 & 4 \\ \hline
    \end{tabular}
    \end{center}
  \caption{Total number of sweeps needed. For the inside partitions  three solve sweeps are performed to created the spikes in the factorization, and four solve sweeps are performed in the solve stage. For the first and last partitions two solve sweeps are performed in the solve stage, and none are required in the factorization stage.}
    \label{tab:sweeps} 
\end{table}

The total number of full sweeps needed for the factorization and solve stages is summarized in Table~\ref{tab:sweeps}.
  We note that in the case where only two partitions are present (i.e. the first and last partition), SPIKE performs the same number of total sweeps than
  a traditional LU factorization and solve would require on solving the original linear system. Since each partition contains half of the elements of the total matrix, 
  a two-partition SPIKE solver that uses one processor/core by partition is expected to run twice faster than a single processor/core LU applied to the whole system
  \cite{Mendiratta:2011}. This is a remarkable result of near perfect parallelism which is often difficult to obtain for complex algorithms due to
  the cost of overhead and additional preprocessing stage.
  This case is known as the SPIKE 2$\times$2 kernel and it will be used as building block in the next sections.

\section{Flexible partitioning scheme for recursive SPIKE} \label{Increased_parallelism_for_Recursive_Spike}

The recursive SPIKE algorithm can only be applied if the number of partitions is a power of two. 
Indeed, the recursive solver repeatedly applies SPIKE to the reduced system, splitting in half the number of partitions with each step.
In previous implementations of recursive SPIKE using OpenMP for shared memory \cite{Mendiratta:2011} or MPI for distributed memory \cite{Polizzi:2011},
the number of threads (resp. MPI processes) was tied to the number of partitions, with one thread (resp. one MPI process) working on each partition.
As a result, the power-of-two restriction for the number of partitions would result in a waste of parallel computing resources. 
For example if 60 cores/processors were available, only 32 cores/processor (the lowest nearest power of two) could be utilized by the standard recursive SPIKE.
The approach discussed in the following waives this restriction by exploiting further the potential for parallelism.
For clarity and without loss of generality (since both MPI and OpenMP SPIKE are possible choices),
the presentation terminology and numerical results are considering a SPIKE OpenMP implementation and the use of threading.

A straightforward method of effectively using additional threads by partition is now proposed. %has been included in this implementation.
If the number of threads is not a power of two, some partitions are given two threads. 
For these partitions the SPIKE 2$\times$2 kernel is used to perform the factorization and solve operations on the associated sub-matrices.
As mentioned in Section~\ref{opt_sweeps}, the SPIKE 2$\times$2 kernel has twice the performance of a single-threaded banded matrix solver. 
Because the factorization and solve operations make up the majority of the computational cost for SPIKE, the 2$\times$2 kernel will provide a significant
speedup for the partitions on which it is used.  

The matrix factorization and solve operations have well known computational costs. 
For banded matrices, the relevant factors are the matrix size and the matrix bandwidth.
The matrix solve operation may also be performed on multiple vectors. 
Load balancing will be achieved by changing the size of each partition so that the computational costs of the large matrix operations on each partition are matched. 
Ultimately this will allow for the definition of optimized ratios between the partition sizes. 

\subsection{Distribution of threads}

This section discusses how threads are allocated to partitions. 
The overall plan is to start by selecting the greatest power of two below the number of available threads to generate the SPIKE partitions,
as is usually the case with recursive SPIKE. 
From there, threads will be added to the middle partitions until we have reached the total number of threads given by the environment. 
Not all partitions will benefit from the addition of threads. 
Specifically, the first and last partitions benefit greatly from exploiting the structure of the LU and UL factorizations respectively, as seen in Section~\ref{opt_sweeps}.
So, conventional LU and UL factorizations are always used for these partitions. 
For all other partitions 2$\times$2 SPIKE may be useful. 

\begin{figure}[htbp] 
\begin{centering}
\includegraphics[keepaspectratio]{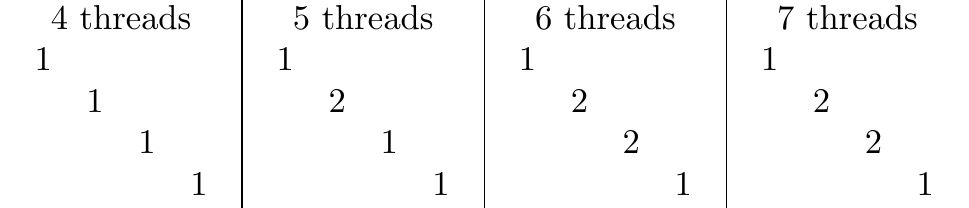} 
\caption{Distribution of 4 to 7 threads using four partitions}
\label{4-7threads}
\end{centering}
\end{figure}

\begin{figure}[htbp] 
\begin{centering}
\includegraphics[keepaspectratio,width=\textwidth]{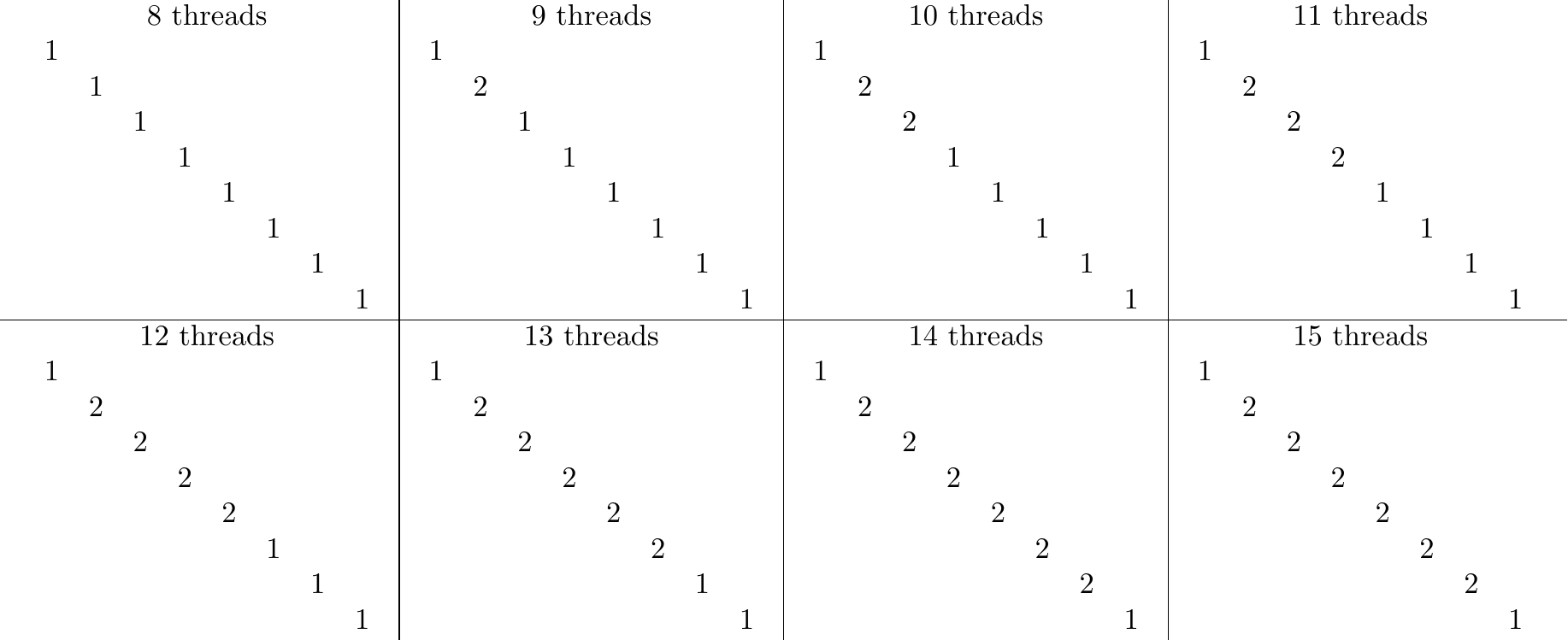} 
\caption{Distribution of 8 to 15 threads using 8 partitions}
\label{8-15threads}
\end{centering}
\end{figure}

Threads are allocated sequentially, starting at the second partition, as shown in Figures~\ref{4-7threads} and~\ref{8-15threads}.
The number one designates a partition which is given a single thread, and the number two designates one given a pair of threads. 
Note that seven threads are distributed as if there were six. 
This is because neither of the remaining single threaded partitions would benefit from using 2$\times$2 SPIKE.
Similarly, in Figure~\ref{8-15threads} one thread is wasted when there are fifteen total threads.
In comparison with the standard recursive SPIKE that allows only one thread per partitions, up to three threads would be wasted
in Figure~\ref{4-7threads} and up to seven in Figure~\ref{8-15threads}.

Formally and in general, we may have $p=2^m$ partitions, for some integer $m$. 
Of these partitions $q$ are given a single thread and $r$ are given two, for a total of $t$ threads.
Finally, the first and last partitions must be given a single thread each. 
Because $m$ is freely selected, any number of threads $t$ may be used with the exception of $2^m-1$ as shown below:

%{{{
\begin{gather}
2 \leq q\leq p;  \quad
0 \leq r \leq p-2, \\ 
p =2^m=q+r,     \\
t =q+2r=(q+r)+r=2^m+r, \\
 2^m \leq t \leq 2^m + p-2 = 2(2^m) - 2,\\
 2^m \leq t \leq 2^{m+1} - 2.
\end{gather}
%}}}

Because the SPIKE algorithm is a domain decomposition algorithm, replacing the \lapack LU solver with a 2$\times$2 SPIKE solver is, algorithmically, trivial.
The derivation of SPIKE given in Section~\ref{basics_section} did not rely on the specifics of the LU factorization, with the exception of a couple of optimizations. 
So, neglecting these optimizations, the 2$\times$2 SPIKE solver may be plugged into place with no changes. 

Of the two main optimizations, only one must require our attention. %significant work to retain. 
The first optimization was used to reduce the number of solve sweeps in the first and last partitions, shown in Section~\ref{opt_sweeps}. 
As stated previously, we simply avoid using the SPIKE 2$\times$2 solver on those partitions, so this is not a problem. 
The more interesting optimization allows for the generation of the $V$ spike beginning with a truncated solve operation, for a total of only one solve-sweep. 
The next section describes how to perform a nearly equivalent optimization, but with the 2$\times$2 SPIKE solver. 

\subsection{Reducing factorization stage sweeps}
In section~\ref{opt_sweeps}, a method of generating the $V$ spikes with just one sweep was shown. 
The essential observation is that the submatrix used to generate $\MBF{V_i}$ is comprised mainly of zeroes, and non-zero elements are restricted to the bottom $\BW$ rows.
As a result, the L-sweep may start at the beginning of the non-zero elements.
This reduces the size of the solve operation from asymptotically equal to the matrix size, to the bandwidth.
As a result it is computationally inexpensive enough to be ignored.

A similar observation can also be applied to the spikes generated with the $2 \times 2$ SPIKE partitions. 
In this case, we will exploit the shape of the $\MBF{B}$ and $\MBF{C}$ matrices to avoid performing solve operations over a large number of zeroes. 
The operations to be performed are:
\begin{equation}
  \MBF{A^{-1}_i}
  \left[ \begin{BMAT}(b,18pt,18pt){c}{bc:c}
\matmid{\MBF{0}}\\ 
 \MBF{\hat{B}_i}
  \end{BMAT} \right]
%\left[
%  \begin{array}{cc}
%    \rvdots \\
 %\MBF{0} \\%& \MBF{0} \\ 
% \MBF{\hat{B}_i} %& \MBF{0}
%\end{array}
%\right]
=
\MBF{V_i}
;
\quad
\MBF{A^{-1}_i}
\left[ \begin{BMAT}(b,18pt,18pt){c}{b:cc}
     \MBF{\hat{C}_i} \\
\matmid{\MBF{0}}
  \end{BMAT} \right]
%\left[
%\begin{array}{cc}
 % %\MBF{0} &
  %\MBF{\hat{C}_i} \\
 % %\MBF{0} &
 % \MBF{0} \\
 % \rvdots
%\end{array}
%\right]
=
\MBF{W_i}.
\end{equation}

$\MBF{A_i}$ is a submatrix of $\MBF{A}$ for which we would like to use $2\times2$ SPIKE\@.
It has a half  bandwidth of $\BW$ and a size of $n_i$.
The relevant equation is:

\begin{equation}
\left[
\begin{array}{cc:cc}
 \multicolumn{2}{c:}{\multirow{2}{*}{$\MBF{A_{i1}}$}} & \MBF{0} & \\
                                             & & \MBF{\hat{B}_{i1}} & \\ \hdashline
 & \MBF{\hat{C}_{i2}} & \multicolumn{2}{c}{\multirow{2}{*}{$\MBF{A_{i2}}$}}  \\
 &  \MBF{0}  & &  
\end{array}
\right]
\left[
\begin{array}{c}
\MBF{\tilde{X}_{i1}} \\
\MBF{X_{i1b}} \\ \hdashline
\MBF{X_{i2t}} \\
\MBF{\tilde{X}_{i2}} 
\end{array}
\right]
=
\left[
\begin{array}{c}
\MBF{\tilde{F}_{i1}} \\
\MBF{F_{i1b}} \\ \hdashline
\MBF{F_{i2t}} \\
\MBF{\tilde{F}_{i2}} 
\end{array}
\right],
%=
%\left[
%\begin{array}{c}
%\matmid{\MBF{F_{i1}}}
%\\ \hdashline
%\matmidt{\MBF{F_{i2}}}
%\\
%\end{array}
%\right],
\end{equation}

where we can extract

\begin{equation}
\MBF{A_{i1}} 
\left[
\begin{array}{c}
\MBF{\tilde{X}_{i1}} \\
\MBF{X_{i1b}} \\
\end{array}
\right]
+
\left[
\begin{array}{c}
\MBF{0} \\
\MBF{\hat{B}_{i1}}
\end{array}
\right]
\MBF{X_{i2t}}
=
\left[
\begin{array}{c}
\MBF{\tilde{F}_{i1}} \\
\MBF{F_{i1b}} \\
\end{array}
\right],
\end{equation}
\begin{equation}
\left[
\begin{array}{c}
\MBF{\tilde{X}_{i1}} \\
\MBF{X_{i1b}} \\
\end{array}
\right]
=
\MBF{A^{-1}_{i1}} 
\left[
\begin{array}{c}
\MBF{\tilde{F}_{i1}} \\
\MBF{F_{i1b}} \\
\end{array}
\right]
-
\left[
\begin{array}{c}
\MBF{0} \\
\MBF{\hat{B}_{i1}}
\end{array}
\right]
\MBF{X_{i2t}}
=
\MBF{U^{-1}_{i1}} 
\left(
\MBF{L^{-1}_{i1}} 
\left[
\begin{array}{c}
\MBF{\tilde{F}_{i1}} \\
\MBF{F_{i1b}} \\
\end{array}
\right]
-
\MBF{L^{-1}_{i,1}} 
\left[
\begin{array}{c}
\MBF{0} \\
\MBF{\hat{B}_{i1}}
\end{array}
\right]
\MBF{X_{i2t}}
\right).
\end{equation}

We may observe that, when solving for $\MBF{V_i}$, $\MBF{F_{i1}}=\MBF{0}$.
The initial L-sweep over this matrix is thus unnecessary.
This saves a solve sweep of height $n_i/2$,

\begin{equation}
\left[
\begin{array}{c}
\MBF{\tilde{V}_{i1}} \\
\MBF{V_{i1b}} \\
\end{array}
\right]
=
\MBF{U^{-1}_{i1}} 
\left(
-
\MBF{L^{-1}_{i1}} 
\left[
\begin{array}{c}
\MBF{0} \\
\MBF{\hat{B}_{i1}}
\end{array}
\right]
\MBF{V_{i2t}}
\right).
\end{equation}

A similar optimization is possible for $\MBF{W_i}$,
This saves another solve sweep of height $n_i/2$, i.e. 

\begin{equation}
\MBF{A_{i2}} 
\left[
\begin{array}{c}
\MBF{X_{i2t}} \\
\MBF{\tilde{X}_{i2}} \\
\end{array}
\right]
+
\left[
\begin{array}{c}
\MBF{\hat{C}_{i2}} \\
\MBF{0} \\
\end{array}
\right]
\MBF{X_{i1t}}
=
\left[
\begin{array}{c}
\MBF{F_{i2b}} \\
\MBF{\tilde{F}_{i2}} \\
\end{array}
\right],
\end{equation}
\begin{equation}
\left[
\begin{array}{c}
\MBF{X_{i2b}} \\
\MBF{\tilde{X}_{i2}} \\
\end{array}
\right]
=
\MBF{A^{-1}_{i2}} 
\left[
\begin{array}{c}
\MBF{F_{i2t}} \\
\MBF{\tilde{F}_{i2}} \\
\end{array}
\right]
-
\left[
\begin{array}{c}
\MBF{\hat{C}_{i2}} \\
\MBF{0} \\
\end{array}
\right]
\MBF{X_{i1t}}
=
\MBF{L^{-1}_{i2}} 
\left(
\MBF{U^{-1}_{i2}} 
\left[
\begin{array}{c}
\MBF{F_{i2t}} \\
\MBF{\tilde{F}_{i2}} \\
\end{array}
\right]
-
\MBF{U^{-1}_{i2}} 
\left[
\begin{array}{c}
\MBF{\hat{C}_{i2}} \\
\MBF{0} \\
\end{array}
\right]
\MBF{X_{i1b}}
\right),
\end{equation}
\begin{equation}
\left[
\begin{array}{c}
\MBF{V_{i2b}} \\
\MBF{\tilde{V}_{i2}} \\
\end{array}
\right]
=
\MBF{L^{-1}_{i2}} 
\left(
-
\MBF{U^{-1}_{i2}} 
\left[
\begin{array}{c}
\MBF{\hat{C}_{i2}} \\
\MBF{0} \\
\end{array}
\right]
\MBF{V_{i1b}}
\right).
\end{equation}

As a result, an amount of work equal to two half-sweeps is saved. 
This means that the total work performed on the SPIKE 2$\times$2 partitions is equal to that of the normal, single threaded partitions. 
In other words, the SPIKE 2$\times$2 kernel may still be used to form the $V$ and $W$ submatrices with three sweeps. 

\subsection{Load balancing scheme} \label{Load_balancinG_scheme_section}

For optimal load balancing, we would like to have each partition take the same amount of time to complete.
This will be approximated by setting equal the sums of the computational costs for the partitions. 
The computational costs considered will be those incurred by the large factorization and solve operations. 

Let us continue using the same banded matrix $\MBF{A}$ with a size of $n \times n$ and a half bandwidth of $\BW$, as well as our collections of vectors $\MBF{F}$ and $\MBF{X}$, sized $n\times n_{rhs}$.
The costs incurred for each partition are summarized in Table~\ref{Computation_Cost_Table}. Note that in the factorization stage, the $\MBF{V}$ and $\MBF{W}$ spikes must be created for the reduced system. 
These require performing solve operations on blocks with widths equal to the lower and upper bandwidths respectively. 
Because the matrix is considered  structurally symmetric (for clarity), these operations are recorded as solve sweeps of width $\BW$. 

%{{{
\begin{table}[htbp] 
\begin{center}
\begin{tabular}{|c||c:c|c|}  \hline
                                 &\multicolumn{3}{l|}{Operation Count}                                         \T{1.7}\\ \cline{2-4}
 \multirow{2}{*}{Partition Type} &\multicolumn{2}{l|}{Factorize Stage}                & Solve Stage            \T{2.0}\\ \cline{2-4}
                                 &\multirow{2}{*}{Factorize} &  Solve Sweeps          & Solve Sweeps            \T{1.8}\\
                                 &                           & (over $\BW$ vectors)   &(over $n_{rhs}$ vectors)\T{1.8}\\ \hline
First \& Last                    &   1                       &  0                     &   2 (LU)               \T{1.8}\\       
Inner Two-Thread                 &   1                       &  3 (SPIKE $2\times2$)&   4 (SPIKE $2\times2$) \T{1.8}\\       
Inner Single-Thread              &   1                       &  3 (LU)              &   4 (LU)               \T{1.8}\\ \hline
\end{tabular}
\end{center}
\caption{Computational cost summary for each partition type.}
\label{Computation_Cost_Table}
\end{table}
%}}}

Table~\ref{Computation_Cost_Table} suggests that one may want to consider three partition sizes, $n_1$, $n_2$, and $n_3$.
Respectively, they are the sizes of the first/last partitions, the middle partitions on which the two threaded SPIKE is used, and the middle partitions which receive the single threaded LU factorization. 
Both types of middle partitions have the same total number of solve sweeps  in each stage. 
The SPIKE $2\times 2$ solver should require half of the computation time used by the standard LU solver. 
So, we may set $n_2 = 2n_3$.
The relationship between $n_{1}$, $n_{2}$ $n_{3}$ can be defined as ratios: $R_{12} = \frac{n_1}{n_2}$ and $R_{13} = \frac{n_1}{n_3}$.

The SPIKE implementation uses a blocked LU factorization and solve, based on the BLAS-3 and \lapack  implementation provided by the system.
Similar to the banded \lapack operations, the factorization has an asymptotic performance of O($n\times \BW^2 $), and the solve has
a performance of O($n\times \BW \times n_{rhs}$). 
These costs can be approximated as $K_1 \times n \times \BW^2$ and $K_2\times n\times \BW \times n_{rhs}$ (using two full sweeps),
and ratio between $K_2$ and $K_1$ may be called $K$.
Because $K$ does not depend on the size of the matrix used, it will become a machine specific tuning constant.
The coefficients $R_{12}$ and $R_{13}$ may be computed by balancing the factorization and solve performance costs
  between the first/last partition and the inner partitions described in Table~\ref{Computation_Cost_Table} as follows:

\begin{gather}
\begin{split}
K_1 n_1 \BW ^2 + K_2 n_1 \BW n_{rhs} = 
K_1n_3\BW ^2+3\frac{K_2}{2}n_3\BW ^2+2K_2n_3\BW n_{rhs} 
\end{split}
\end{gather}
\begin{gather}
\begin{split} 
K_1 n_1 \BW + K_2 n_1  n_{rhs} = 
K_1n_3\BW +(3/2)K_2n_3\BW +2K_2n_3n_{rhs}
\end{split}.
\end{gather}
Now it is possible to obtain $R_{13}$ in terms of $K$, $n_{rhs}$, and $\BW$:
\begin{gather}
K_1 n_1 \BW ^2 + K_2 n_1 \BW n_{rhs} = 
K_1n_3\BW ^2+(3/2)K_2n_3\BW ^2+2K_2n_3\BW n_{rhs} 
\\
\begin{split} \label{R13_calculation}
R_{13} &= \frac{n_1}{n_3} =
\frac
{K_1\BW +(3/2)K_2n_3\BW +2K_2n_{rhs}}
{K_1 \BW + K_2 n_{rhs}}
\\&=
\frac
{1}
{1 + (K_2/K_1)(n_{rhs}/\BW)}
+
\frac
{3/2+2n_{rhs}/\BW}
{K_1/K_2  + n_{rhs}/\BW}
\\&=
\frac
{1}
{1 + (K)(n_{rhs}/\BW)}
+
\frac
{3/2+2n_{rhs}/\BW}
{1/K  + n_{rhs}/\BW}.
\end{split}
\end{gather}

For $R_{12}$ we have:
\begin{gather}
n_2 = 2 n_3, \\
R_{12} = \frac{1}{2}R_{13} = \label{R12_calculation}
\frac
{1}
{2 + 2(K)(n_{rhs}/\BW)}
+
\frac
{3/4+n_{rhs}/\BW}
{1/K  + n_{rhs}/\BW}.
\end{gather}

The constant $K$ depends on the system hardware and the underlying \lapack and BLAS implementations.
Due to the myriad of existing hardware and software, it is unlikely that an universally good value for $K$ exists. 
However, for a given machine $K$ may be easily found by performing a matrix factorization and solve on a matrix and set of vectors for which $n_{rhs}=\BW$. 
Using the same approximations as above,

\begin{gather}
\mbox{factorization time} = K_1 \times n \times \BW^2 \label{ftime}, \\ 
\mbox{solve time}= K_2\times n\times \BW \times n_{rhs} \label{stime}, \\
K=\frac{K_2}{K_1}
 =\frac{\mbox{solve time}}{n \times \BW \times n_{rhs}} \times \frac{n \times \BW^2}{\mbox{factorization time}} \\
 =\frac{\mbox{solve time}}{\mbox{factorization time}}.
\end{gather}

This calculation requires that the matrix used is large enough for the asymptotic computational costs to dominate. 
The implementation of SPIKE discussed here contains the ability to include a value for $K$ as an input parameter. 
Because $K$ is constant for a given machine and BLAS/LAPACK implementation, it could be computed once and for all
  after installation of the SPIKE software package.

The other variable to consider when determining $R_{12}$ and $R_{13}$ is $n_{rhs}/\BW$.
In general, if this value is known before the DS factorization is performed, $R_{12}$ and $R_{13}$ may be calculated. 
If the value is not known, the problem might be characterized as similar to one of two limiting cases, $n_{rhs}/\BW \rightarrow 0$ and $n_{rhs}/\BW \rightarrow \infty$.

In the first case, the matrix bandwidth is much greater than the number of vectors in the solution.
Intuitively, this indicates that the factorization stage will dominate the computational cost.
In this case, we obtain:
\begin{gather}\label{ratiof}
\mbox{lim}_{n_{rhs}/\BW} \rightarrow 0
\hspace{10pt}
R_{12} = (1/2)+(3/4)K.
\end{gather}
This can be seen simply by plugging the value $n_{rhs}/\BW=0$  into equation  (\ref{R12_calculation}) for $R_{12}$.

In the second case, where the number of solution vectors is much greater than the matrix bandwidth, the solve stage dominates. 
For this type of problem, we obtain constant ratios that are independent of the value of $K$ i.e.
\begin{gather}\label{ratios}
\mbox{lim}_{n_{rhs}/\BW} \rightarrow \infty
\hspace{10pt}
R_{12} = 
\frac{1}{2 + 2(K) (n_{rhs}/\BW)}
+ 
\frac{1 + n_{rhs}/\BW}{1/K + n_{rhs}/\BW}
=1, \ \mbox{and} \quad R_{13}=2.
\end{gather}

Once the ratios between partition sizes have been decided upon, sizing the partitions is simple. 
The main requirement is that the partition sizes must sum to the size of $\MBF{A}$. 
Assuming there are $x=r-2$ partition of size $n_2$, $y\equiv q$ of size $n_3$, and the first and last partitions, each of which is size $n_1$. 
Overall, this gives the following constraints, which can be trivially solved for the size of each type of partition:

\begin{equation}
n = 2n_1 + xn_2 + yn_3
= 2n_1 + \frac{xn_1}{R_{12}} + \frac{yn_1}{R_{13}},
\end{equation}
\begin{gather}
\frac{nR_{12}R_{13}}{2R_{12}R_{13}+xR_{13}+yR_{12}}= n_1,  \\
\frac{nR_{13}}{2R_{12}R_{13}+xR_{13}+yR_{12}}= n_2,  \\
\frac{nR_{12}}{2R_{12}R_{13}+xR_{13}+yR_{12}}= n_3.
\end{gather}

This concludes the description of the increased parallelism scheme for recursive SPIKE\@.
In summary, this scheme allows the use of almost any number of threads, without dramatically modifying the recursive SPIKE algorithm.
Overall computational time is decreased by carefully sizing the partitions into which the matrix $\MBF{A}$ is distributed. 
The information required for the sizing process has been separated into hardware/library-dependent factors and problem-dependent ones.
Finally, the sizing task is simple enough that it may be performed automatically, and the SPIKE OpenMP library package \cite{spike} includes utility routines to do so. 

\subsection{Performance measurements}\label{sec:machine}
To show the effects of the previously described enhancements, a number of measurements were taken on a large shared memory machine. 
The first set of measurements explore the partition sizing method, as described in the previous section. 
The second set of measurements shows the overall performance and scalability of the algorithm.
The hardware and software used for these experiments is as follows: 
\begin{itemize}
\item{8$\times$Intel\textregistered\ Xeon\textregistered\ E7-8870: 10 cores @ 2.40 GHz with 30MB cache}
\item{Intel\textregistered\ Fortran 16.0.1}
\item{Intel\textregistered\ MKL 11.3.1}
\end{itemize}

The E7-8870 also exploits the `hyperthreading' simultaneous multithreading strategy.
Hyperthreading is generally considered to be detrimental for dense numerical linear algebra.
In most cases, for these experiments hyperthreads have been avoided using the following environment variable:
\begin{itemize}
\item{KMP\_AFFINITY=granularity=fine,compact,1,0}
\end{itemize}
The KMP affinity interface is a feature of the Intel implementation of
OpenMP.\footnote{By default, the pair of hyperthreads run by a given CPU core are considered to be hierarchy very close to one-another. 
The `compact' command instructs the OpenMP runtime to pack threads as closely as possible.
The `1,0' command shifts the core hierarchy, so that the pair of hyperthreads on a given core are considered very far away from one another, while the cores inside a given CPU package are considered nearest neighbors. 
By using this strategy and employing less than eighty threads, a pair of hyperthreads which share a core are never considered close enough to employ both simultaneously.} 

Finally, SPIKE is also making extensive use of LAPACK/BLAS3, so any improvements in the kernel library (e.g. Intel MKL) would be as well beneficial to SPIKE 
and it would not change the relative scalability and speed-up performances  between SPIKE-OpenMP and MKL that are presented here.

\subsubsection{Partition ratio accuracy}

\begin{figure}[htbp] 
\includegraphics[keepaspectratio,width=\textwidth]{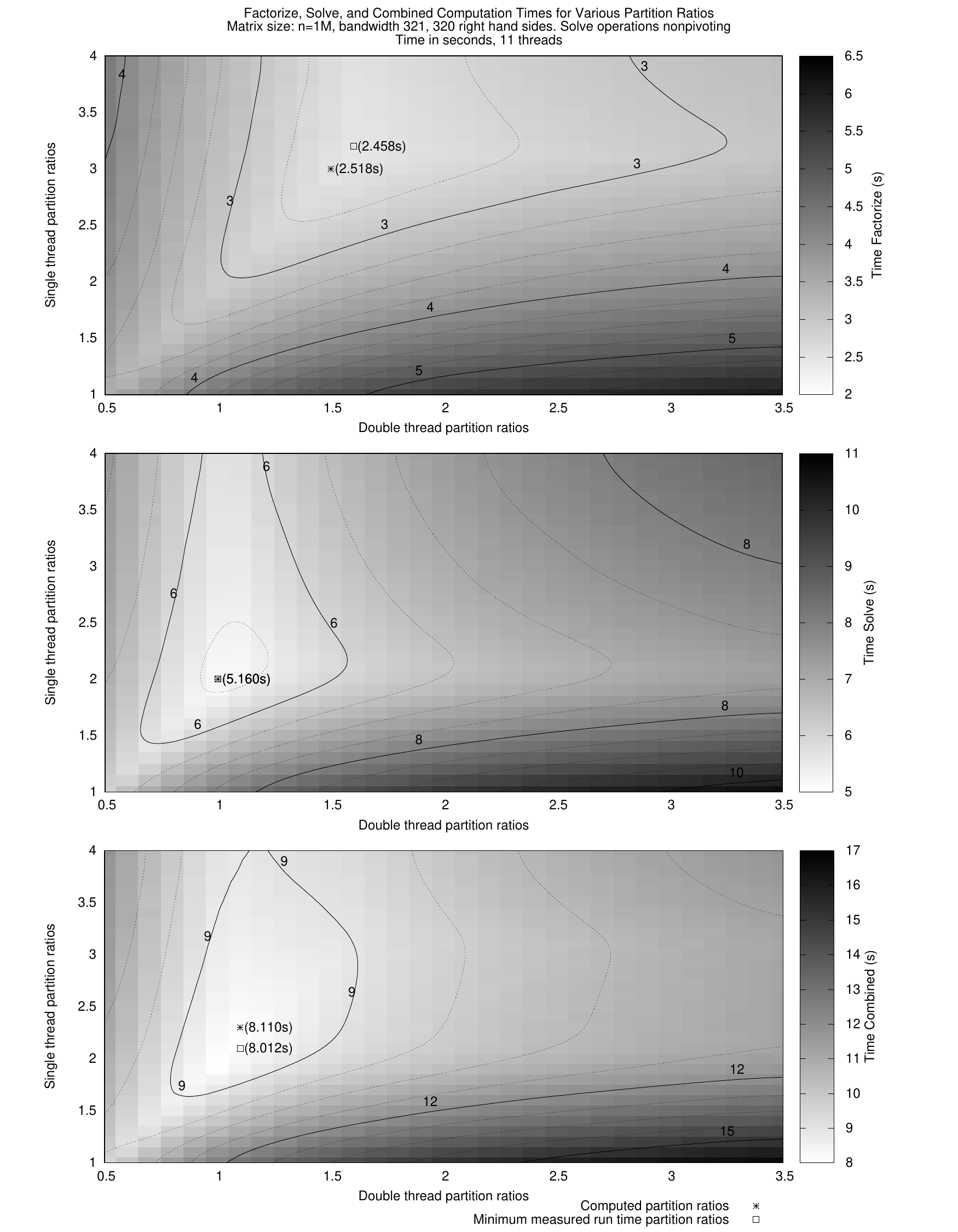} 
\caption{Partition ratio `heatmaps' for 320 right hand sides}
\label{heatmapNonpivot320}
\end{figure}

\begin{figure}[htbp] 
\includegraphics[keepaspectratio,width=\textwidth]{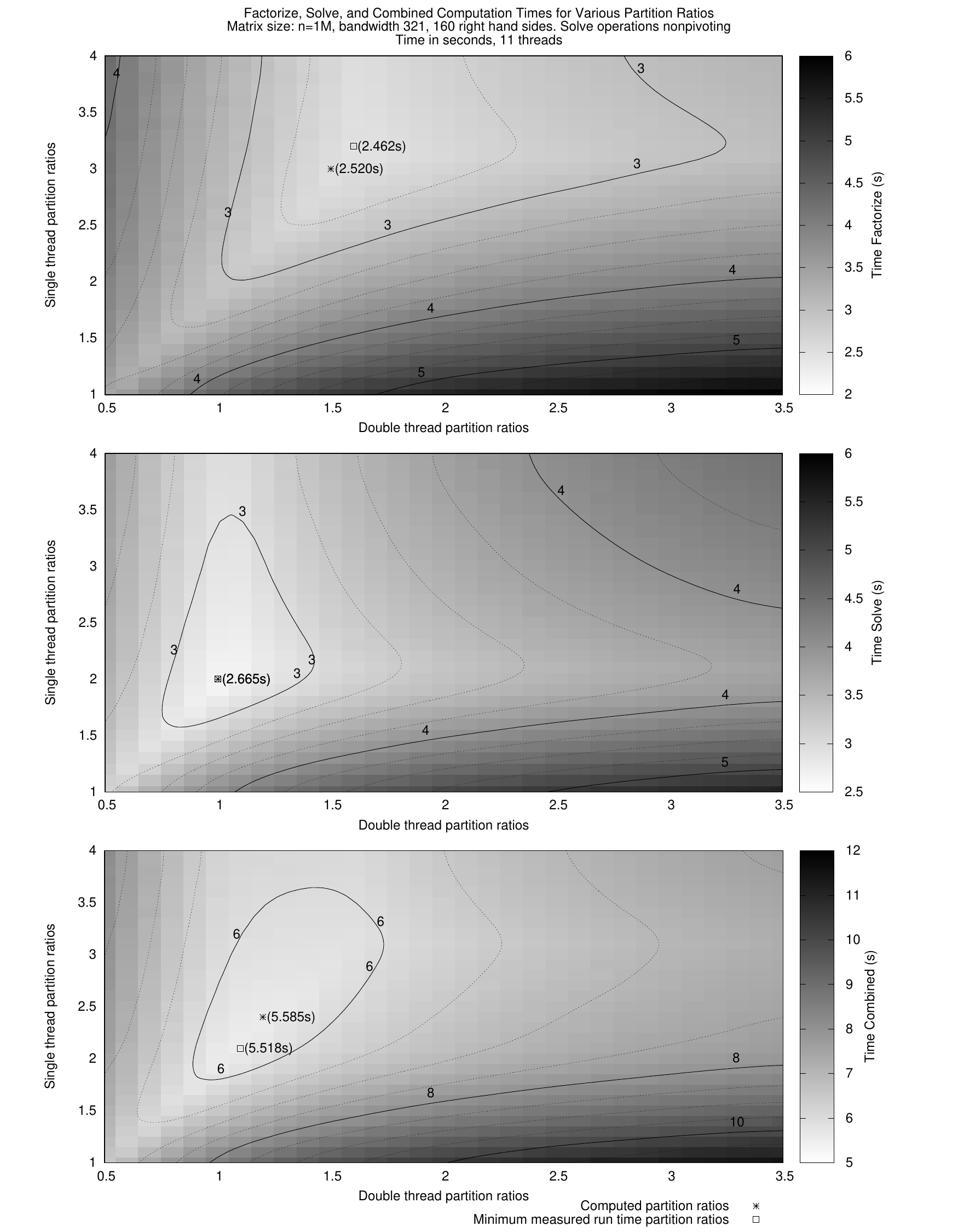} 
\caption{Partition ratio `heatmaps' for 160 right hand sides}
\label{heatmapNonpivot160}
\end{figure}

\begin{figure}[htbp] 
\includegraphics[keepaspectratio,width=\textwidth]{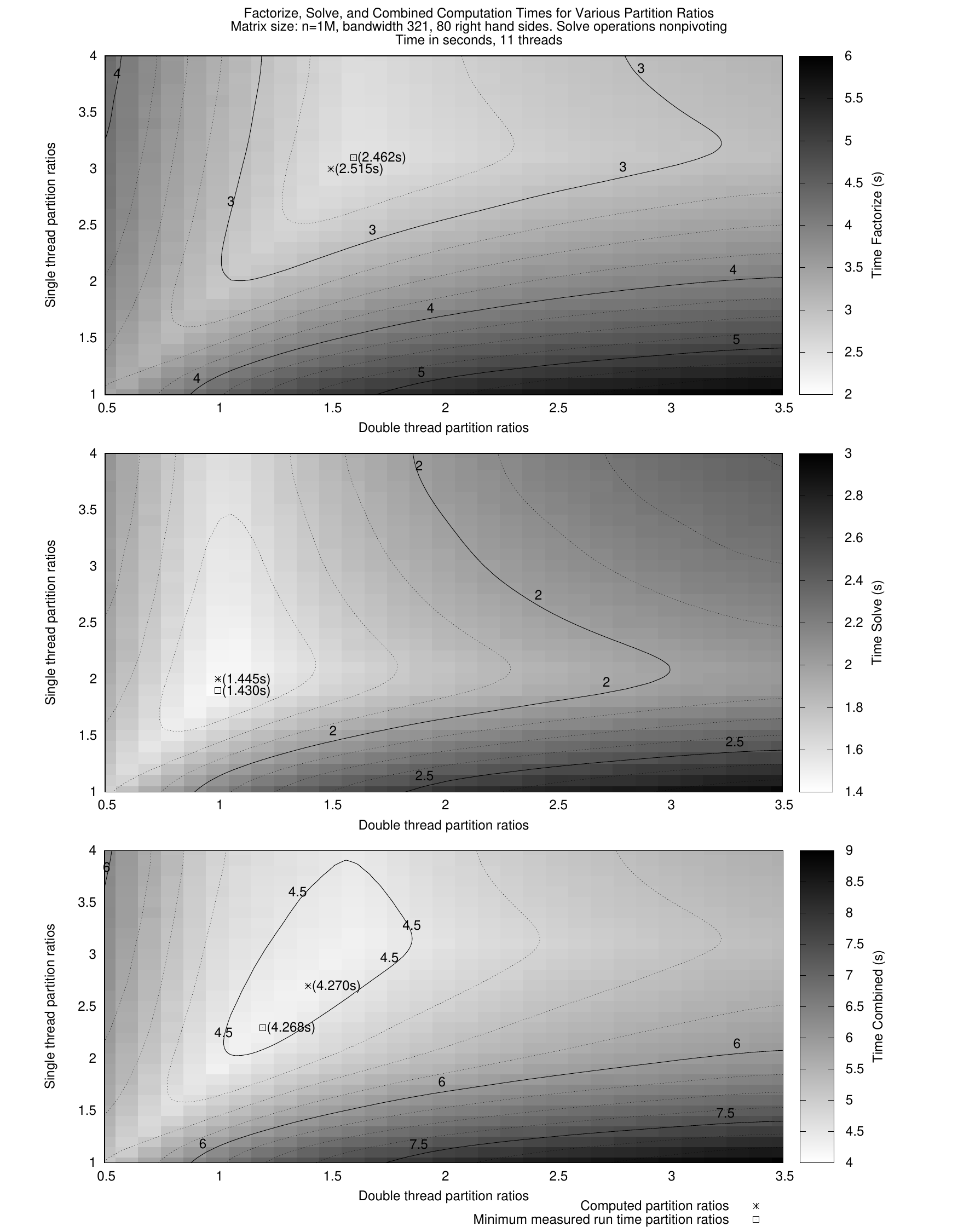} 
\caption{Partition ratio `heatmaps' for 80 right hand sides}
\label{heatmapNonpivot80}
\end{figure}

In Section~\ref{Load_balancinG_scheme_section} equations to determine the appropriate sizes of the various submatrices used in the domain decomposition are derived.
To measure the accuracy of this technique, an exploration of many possible partition size ratios was performed in 
Figures~\ref{heatmapNonpivot320} through~\ref{heatmapNonpivot80}.
For these measurements, the matrix size $n$ and bandwidth $b$ remain constant (resp. $n=10^6$ and $b=321$ with $\BW=160$), while the number of solution vectors changes from $n_{rhs}=320$ in Figure~\ref{heatmapNonpivot320} to
$n_{rhs}=160$ in Figure~\ref{heatmapNonpivot160}, and then $n_{rhs}=80$ in Figure~\ref{heatmapNonpivot80}.  %is varied. 
In these figures, the X and Y axes correspond to the ratios $R_{12}$ and $R_{13}$, as defined in Section~(\ref{Load_balancinG_scheme_section}).
By keeping the bandwidth constant and varying the number of solution vectors, the effect of these ratios can be observed.
Each figure has a map for the cost, in time, of the factorization and solve stages, as well at the overall computation time.
In addition, the best measured runs as well as the location of the pre-calculated values of the best partition size ratios, have been marked along with their times.
The pre-computed values for the factorization and solve stages use the most favorable ratios derived in (\ref{ratiof}) and (\ref{ratios}), respectively.
  The pre-computed value for the combined factorization/solve measurement is obtained
  using the ``compromise ratios'' given in equations~(\ref{R13_calculation}) and~(\ref{R12_calculation}).
Because the matrix does not change from one run to the next, the factorization stage is identical for each run. 
As such, the first map in each figure is largely identical, with some small variation due to noise. 
The excellent agreement between the results indicates 
  that $K$ the machine specific tuning constant,
  is  accurately computed.
The method of determining the solve stage favoring partition ratios is even more reliable than the factorization stage. 
Indeed, for Figures~\ref{heatmapNonpivot320} and~\ref{heatmapNonpivot160} the measured and calculated values are identical. 
This is likely because the solve stage partition ratio formula can be simplified to a pair of constant numbers, so whatever imprecision was introduced
in the discovery of $K$ is no longer present.
%%%%%%%%%%%%%%%%%%%%%%%%%%%%%%%%%%%%%%%%%

Finally, a band of good performance can be visually observed starting at the origin and continuing along the path of $2R_{12} = R_{13}$.
Within those areas, the primary concern is that the computation times produced by using the calculated partition ratios are not too far from the optimal measured ones.
The percentage improvement from using the measured optimal, rather than calculated, partition ratios is summarized in Table~\ref{Computation_Percentage_Table}
for $n_{rhs}=40$ to $n_{rhs}=320$.
In general the gains of the measured optimal partition ratios over the computed ones are in the low single-digit percentages.

\begin{table}[htbp] 
\begin{center}
\begin{tabular}{|c||c|c|c|c|}  \hline
Solution Vectors &   40   &   80   &   160  &   320  \\ \hline
Factorize        & 2.44\% & 2.15\% & 2.36\% & 2.44\% \\ 
Solve            & 1.43\% & 1.05\% &   0    &   0    \\ 
Combined         & 1.22\% & 0.04\% & 1.21\% & 1.22\% \\ \hline
\end{tabular}         
\end{center}
\caption{Performance gain from using best measured partition ratios $\left|\frac{t_{calculated}}{t_{measured}}-1\right|$}
\label{Computation_Percentage_Table}
\end{table}

\subsubsection{Scalability and performance comparisons}\label{Non_Pivoting_scaling}

\begin{figure}[htbp] 
\includegraphics[keepaspectratio,width=\textwidth]{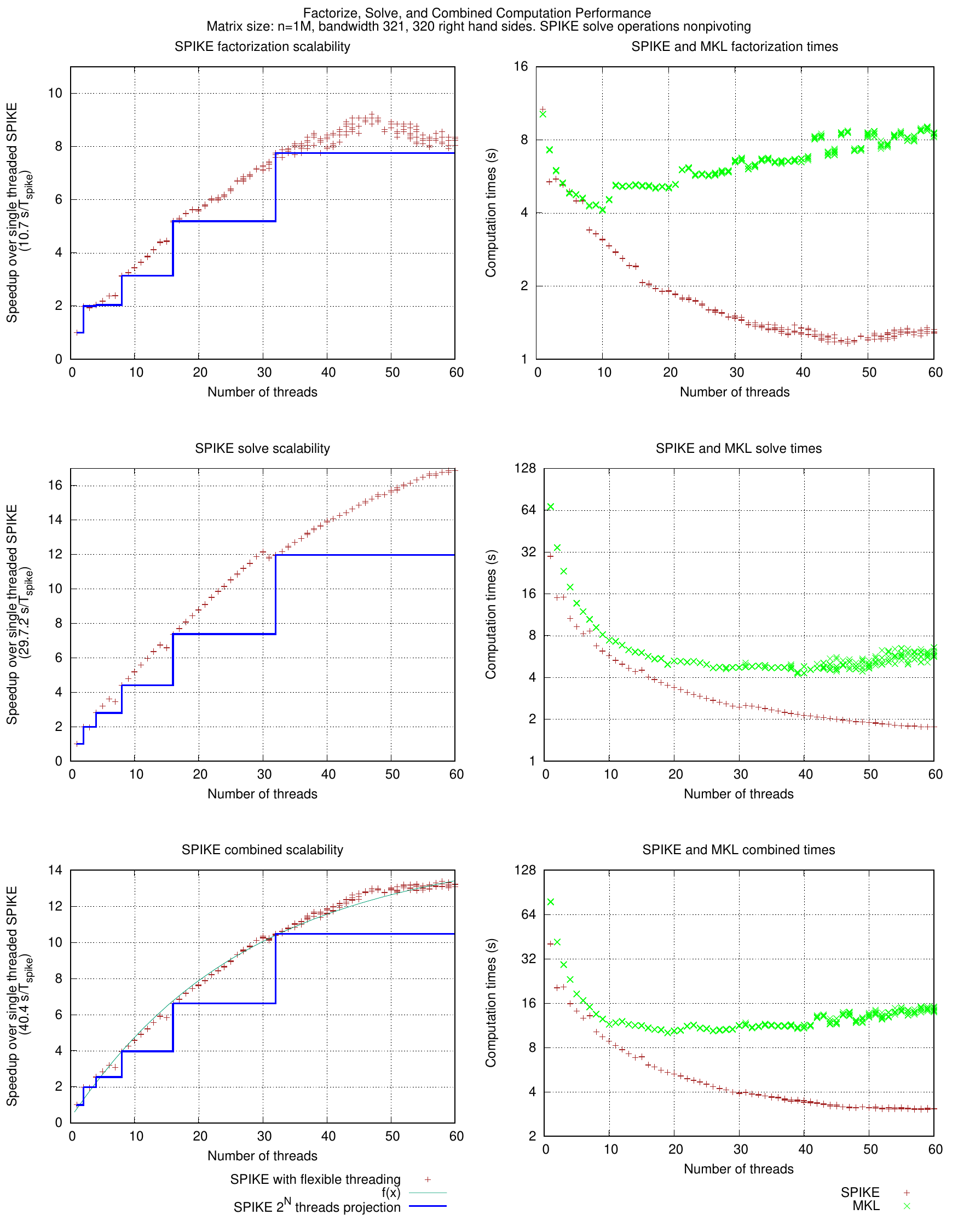} 
\caption{Scalability and computation time for 320 right hand sides}
\label{scalabilityNonpivot320}
\end{figure}

\begin{figure}[htbp] 
\includegraphics[keepaspectratio,width=\textwidth]{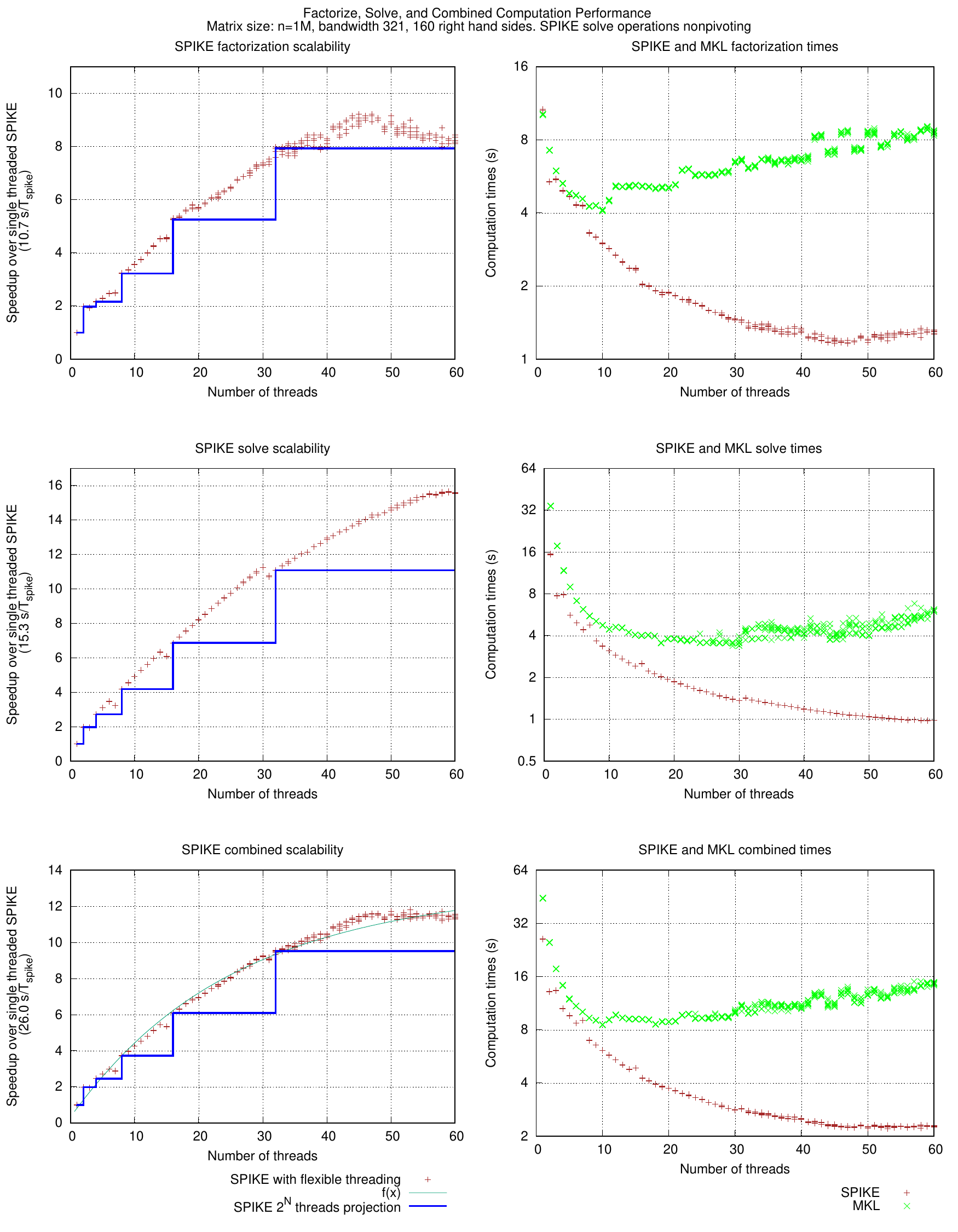} 
\caption{Scalability and computation time for 160 right hand sides}
\label{scalabilityNonpivot160}
\end{figure}

\begin{figure}[htbp] 
\includegraphics[keepaspectratio,width=\textwidth]{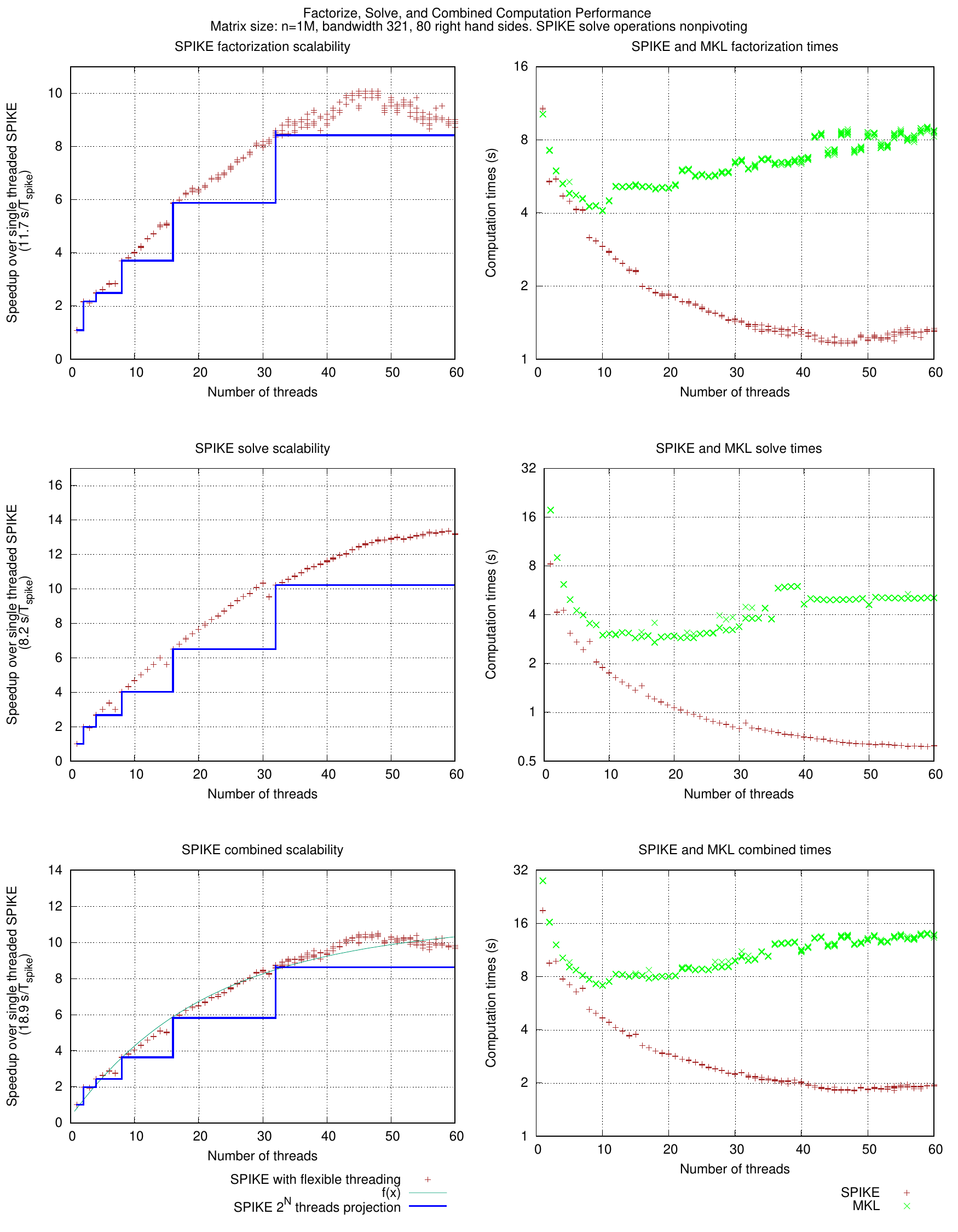} 
\caption{Scalability and computation time for 80 right hand sides}
\label{scalabilityNonpivot80}
\end{figure}

We propose to  observe some aspects of the overall performance of the new implementation of recursive SPIKE.
Figures~\ref{scalabilityNonpivot320},~\ref{scalabilityNonpivot160}, and~\ref{scalabilityNonpivot80} contain two sets of measurements. 
On the left, we see the scalability of SPIKE. 
On the right, we see absolute time measurements, as well as a comparison to MKL (Note that the time axes in these measurements are on a logarithmic scale).
All measurements for SPIKE (including factorization, solve and combined stages) were taken using the calculated partition ratios given in ~(\ref{R13_calculation}) and~(\ref{R12_calculation})
and summarized in Table~\ref{tab:ratios}.

\begin{table}[htbp]
  \begin{center}
\begin{tabular}{|c||c|c|c|}  \hline
Solution Vectors  & 80  & 160 & 320 \\ \hline
R13               & 2.7 & 2.4 & 2.3 \\ 
R12               & 1.35 & 1.2 & 1.15 \\ \hline
\end{tabular}         
  \end{center}
  \caption{Partition ratios used for Figures~\ref{scalabilityNonpivot320},~\ref{scalabilityNonpivot160}, and~\ref{scalabilityNonpivot80}.}
\label{tab:ratios}
\end{table}

Scalability is measured relative to the computation time of the single-threaded non-pivoting solver used on the individual partitions. 
Overall, scaling for the combined factorization/solve stages, continues quite well until around 45 cores are used. 
After that point, the results stall and would eventually degrade in performances. 
We note that the scalability breaking point could go well beyond the 45 cores while considering larger matrices.
The trade-off used to determine the partition ratios can be seen by comparing the scaling of each set of benchmarks. 
As the number of solution vectors decreases, the partition size ratios move to favor the factorization stage of the computation. 
This can be observed in the increased scaling of the factorization stage, and the decrease in the solve stage scaling.
We note that 
  the optimal ratios for the factorization stage given in (\ref{ratiof}) are equal to  $R_{13}=3$ and $R_{12}=1.5$ for
 the measured value of $K$ on our software/hardware set-up. 
The ratios
provided in Table~\ref{tab:ratios} will progressively reach these values with the number of solution vectors decreasing.
In turn, the optimal ratio for the solve stage (\ref{ratios}) give  the values $R_{13}=2$ and $R_{12}=1$, which are close to the values reported in Table~\ref{tab:ratios}
  with large number of right hand sides. Overall for these particular numerical experiments, the solve stage has noticeably superior scalability to the factorization stage.

The scalability measurements also show the benefit of the flexible threading scheme.
This is one of the most important results presented here, since the standard recursive SPIKE scheme is limited by the use of power of two
  number of threads. The line labeled `SPIKE $2^N$ threads projection' shown the effects of limiting the number of threads used
to powers of two by extending the performance measured at these points.
Naturally, the performance gap is most dramatic soon before the number of threads is increased to the next power of two. 
For example, looking at Figure~\ref{scalabilityNonpivot160}, at 30
threads the overall computation scaling increases from roughly 6$\times$ to roughly 9$\times$, as a result of the increased overall utilization of resources.

Finally, overall computation time is generally superior to MKL.
We note that the two solvers are close in time until 10 threads are reached, at which point SPIKE begins pulling away. 
This is particularly apparent in the factorization stage.  In contrast to the SPIKE $DS$ factorization, parallelism performance for
  the inherently recursive serial $LU$ approach used by MKL  mainly relies on BLAS which quickly reaches its limits.
On the other hand, MKL parallelizes well over solution vectors, and so when their number increases, MKL remain moderately
closer in performance to SPIKE.
We note that the base solver used for SPIKE provides performance advantage, as it is non-pivoting. In order to minimize the effects
  of pivoting for MKL, all the test matrices in the numerical experiments were chosen diagonally dominant
(both solvers producing relative residuals of $10^{-13}$ or below). However, SPIKE recursive is applicable
  to non-diagonally dominant systems as well. In most cases, a zero-pivot may never been found 
  even for matrices with large condition numbers. The latter, however, could affect the relative residual and a SPIKE pivoting strategy
  will be presented in Section~\ref{Pivoting_Spike_Section} to address this issue.

\subsubsection{Comments on hardware specific configuration}\label{hardware}

In all cases, the problem is configured on the master thread; that is, core 0 of CPU package 0. 
Memory is allocated in this thread.
This is representative of the expected use case for this code. 
Our intention is to create a black-box shared memory implementation of SPIKE.
It is unreasonable to expect a user to reconfigure their code -- formulate the creation of their matrices in parallel -- simply to
replace the matrix factorize and solve operations. 
However, this does cause what appear to be some non-uniform memory access (NUMA) issues. 
In particular, it seems that the CPU package 0 and 1 have faster access to memory allocated by cores on CPU 0. 
These issues were not apparent in the previous section; the E7-8870 is a 10-core CPU, and so with 11 threads a compact core allocation
method results in the cores being limited to CPU packages 0 and 1.

To minimize these issues, the OpenMP threads were explicitly mapped to the hardware cores. 
A modified `scatter' technique has been used, to maximize locality and cache utilization. 
The number of cores used per CPU is determined by  dividing the number of requested threads by the number of CPU packages ($n_{threads}/8$),
with the remainder simply allocated sequentially starting at CPU 0. 
First, the threads working on the first and last partitions were mapped to the cores 0 and 1 of CPU 0. 
Because these partitions have the least work per element their performance becomes memory bound most rapidly, so locating them on the CPU with
the best NUMA access improves performance.
Next, threads are mapped to cores sequentially using the threads per package rule. 
So, for example, with 16 threads, and thus 16 partitions, we would have partitions 0 and 15 on package 1, 1 and 2 on package 1, 3 and 4 on package 2,
and so on. 
This maximizes the availability of CPU cache (particularly important on a system with a relatively large 30MB of level 2 cache per CPU package) and
NUMA friendliness while minimizing the amount of intra-package communication that must occur when information is passed from one partition to the other.

%%%%%%%%%%%%%%%%%%%%%%%%%%%%%%%%%%%%%%%%%%%%%%%%%%%%%%%%%%%%%%%%%%%%%%%%%%%%%%%%%%%%%%%%%%%%%%%%
%%%%%%%%%%%%%%%%%%%%%%%%%%%%%%%%%%%%%%%%%%%%%%%%%%%%%%%%%%%%%%%%%%%%%%%%%%%%%%%%%%%%%%%%%%%%%%%%

\section{Transpose solve option for recursive SPIKE}

A transpose solve option is a standard feature for \lapack subroutines.
This option allows transpose problems to be solved without explicitly transposing the matrix in memory.
Transpose solve retrieves $\MBF{X}$ for the following problem:
\[\MBF{A}^T\MBF{X}=\MBF{F},\]
where $\MBF{A}$, $\MBF{X}$ and $\MBF{F}$ are defined as in the previous sections: An $n\times n$ banded matrix with half-bandwidth $k$, and two $n\times n_{rhs}$ collections of vectors, respectively.

Similarly to the standard \lapack solver, the transpose solve option reuses the factorization from the non-transpose case. 
That is, once a matrix has been factorized it may be used for either transpose or non-transpose solve operations.
Because the factorization stage has the potential to be much more time-consuming than the solve stage, this feature can result in great time savings. 
For SPIKE, this means we reuse the $\MBF{D}$ and $\MBF{S}$ matrices and the reduced system from the previous section.
The transpose problem may be written as follows:

\begin{gather}
\MBF{A}^T\MBF{X}=\MBF{(DS)}^T\MBF{X}=\MBF{S}^T\MBF{D}^T\MBF{X}=\MBF{F}, \\
\MBF{S}^T\MBF{Y}=\MBF{F}, \\
\MBF{D}^T\MBF{X}=\MBF{Y}.
\end{gather}

This presents two sub-problems. 
As in the non-transpose case, partitions of the $\MBF{D}$ matrix are uncoupled, and so the $\MBF{D^T}$ stage can be parallelized in a familiar, straightforward manner.
For the $\MBF{S}^T$ matrix a new algorithm will need to be designed because this matrix is structurally different from the $\MBF{S}$ matrix.
In particular, a transpose version of the recursive reduced system solver is required. 
Ultimately near performance parity with the non-transpose solver will be achieved by matching the count of these operations. 
This will guide the development of the algorithm. 

\subsection{Transpose S stage} \label{Transpose_S_stage}

\begin{figure}[htbp] 
\begin{centering}
\includegraphics[keepaspectratio]{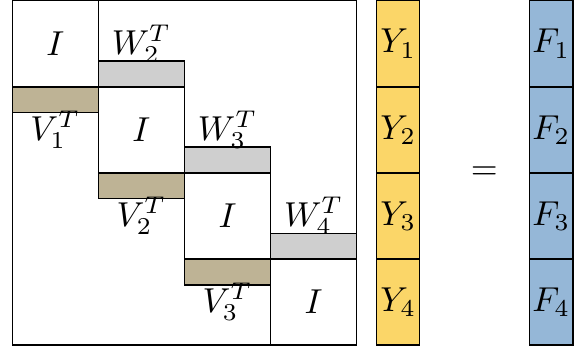} 
\caption{Four partition transpose S-matrix}
\label{S_transpose_1}
\end{centering}
\end{figure}

The first sub problem to solve is $\MBF{S^T}\MBF{Y}=\MBF{F}$.
This problem can be visualized using the four-partition example in Figure~\ref{S_transpose_1}. 
A reduced system can be extracted from this matrix, by exploiting the fact that many of the elements of the $\MBF{Y}$ vector are not affected by the solve operation, and therefore are simply equal to the corresponding elements of $\MBF{F}$. 
This can be seen if the $\MBF{V_i^T}$ and $\MBF{W_i^T}$ spikes, and the $\MBF{Y_i}$ and $\MBF{F_i}$ vectors are partitioned in the following manner:

\begin{gather}
\MBF{V_i}^T=
\left[
\begin{array}{c c c}
\MBF{V_{it}^T}, & \MBF{\tilde{V}_i}^T, & \MBF{V_{ib}^T}
\end{array}
\right]; \quad
\MBF{W_i^T}=
\left[
\begin{array}{c c c}
\MBF{W_{it}^T}, & \MBF{\tilde{W}_i}^T, & \MBF{W_{ib}^T}
\end{array}
\right],
\end{gather}

\begin{gather}
\MBF{Y_i}=
\left[
\begin{array}{c c c}
\MBF{Y_{it}^T}, & \MBF{\tilde{Y}_i}^T, & \MBF{Y_{ib}^T}
\end{array}
\right]^T; \quad
\MBF{F_i}=
\left[
\begin{array}{c c c}
\MBF{F_{it}^T}, & \MBF{\tilde{F}_i}^T, & \MBF{F_{ib}^T}
\end{array}
\right]^T,
\end{gather}

\begin{figure}[htbp] 
\begin{centering}
\includegraphics[keepaspectratio,width=\textwidth]{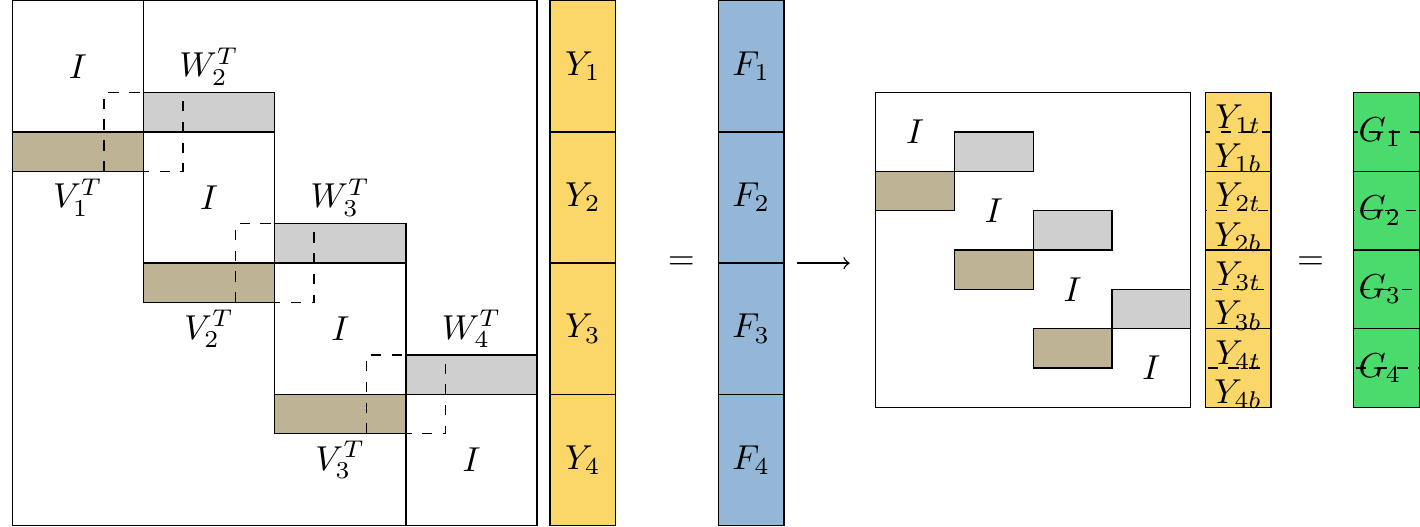} 
\caption{Reduced transpose system extraction for four partitions}
\label{reduced_transpose_1}
\end{centering}
\end{figure}

When viewing a given horizontal slice of the $\MBF{S^T}$ matrix, shown in Figure~\ref{S_transpose_1}, it is visually clear that $\MBF{\tilde{F}_i} = \MBF{\tilde{Y}_i}$.
Indeed, we obtain:

\begin{gather}
\label{Sslice}
\MBF{F_i} = 
\left[
\begin{array}{c}
\MBF{F_{it}} \\ 
\MBF{\tilde{F}_i} \\ 
\MBF{F_{ib}}
\end{array}
\right]
=
\left[
\begin{array}{c}
\MBF{V_{i-1}^T} \\
\MBF{0}\\ 
\MBF{0}
\end{array}
\right] 
\left[
\begin{array}{c}
\MBF{Y_{i-1t}} \\ 
\MBF{\tilde{Y}_{i-1}} \\ 
\MBF{Y_{i-1b}}
\end{array}
\right]
+
\left[
\begin{array}{c c c}
\MBF{I} & \MBF{0} & \MBF{0}\\ 
\MBF{0} & \MBF{I} & \MBF{0}\\ 
\MBF{0} & \MBF{0} & \MBF{I}
\end{array}
\right] 
\left[
  \begin{array}{c}
    \MBF{Y_{it}} \\ 
    \MBF{\tilde{Y}_{i}} \\ 
    \MBF{Y_{ib}}
  \end{array}
\right]
+
\left[
  \begin{array}{c}
    \MBF{0}\\ 
    \MBF{0}\\ 
    \MBF{W_{i+1}^T}
  \end{array}
\right] 
\left[
  \begin{array}{c}
    \MBF{Y_{i+1t}} \\ 
    \MBF{\tilde{Y}_{i+1}} \\ 
    \MBF{Y_{i+1b}}
  \end{array}
\right].
\end{gather}
If $\MBF{Y_{it}}$ and $\MBF{Y_{ib}}$ are given a height of $k$ rows each, and $\MBF{\tilde{Y}_i}$ is given the remaining elements, this equation can be rewritten as follows:
\begin{gather}
\begin{split}
\MBF{F_{it}} &= 
\MBF{Y_{it}} + 
\MBF{V_{i-1}^T}
\left[
  \begin{array}{c}
    \MBF{Y_{i-1t}} \\ 
    \MBF{\tilde{Y}_{i-1}} \\ 
    \MBF{Y_{i-1b}}
  \end{array}
\right]
=
\MBF{Y_{it}} +
\MBF{V_{i-1}^T}
\left[
  \begin{array}{c}
    \MBF{0} \\
    \MBF{\tilde{Y}_{i-1}} \\ 
    \MBF{0} \\
  \end{array}
\right]
+
\MBF{V_{i-1}^T}
\left[
  \begin{array}{c}
    \MBF{Y_{i-1t}} \\ 
    \MBF{0}\\
    \MBF{0}
  \end{array}
\right]
+
\MBF{V_{i-1}^T}
  \left[
    \begin{array}{c}
    \MBF{0}\\
    \MBF{0}\\
    \MBF{Y_{i-1b}}
  \end{array}
\right],
\\
\MBF{\tilde{F}_i}
&=
\MBF{\tilde{Y}_i},
\\
\MBF{F_{ib}} &= 
\MBF{Y_{ib}} + 
\MBF{W_{i+1}^T}
\left[
  \begin{array}{c}
    \MBF{Y_{i+1t}} \\ 
    \MBF{\tilde{Y}_{i+1}} \\ 
    \MBF{Y_{i+1b}}
  \end{array}
\right]
=
\MBF{Y_{ib}} +
\MBF{W_{i+1}^T}
\left[
  \begin{array}{c}
    \MBF{0} \\
    \MBF{\tilde{Y}_{i+1}} \\ 
    \MBF{0} \\
  \end{array}
\right]
+
\MBF{W_{i+1}^T}
\left[
  \begin{array}{c}
    \MBF{Y_{i+1t}} \\ 
    \MBF{0}\\
    \MBF{0}
  \end{array}
\right]
+
\MBF{W_{i+1}^T}
\left[
  \begin{array}{c}
    \MBF{0}\\
    \MBF{0}\\
    \MBF{Y_{i+1b}}
  \end{array}
\right].
\end{split}
\end{gather}
The solve for $Y_{it}$ and $Y_{ib}$ must now  
  be modified to adjust for the presence of the known values in $\MBF{\tilde{Y}_i}$.
It is then possible to extract a reduced system as depicted in Figure~\ref{reduced_transpose_1}, and where the
modified right-hand side $\MBF{G_i}$ is given by:
\begin{gather} \label{gform1}
\begin{split}
i>1, \quad \MBF{G_{it}} =
\MBF{F_{it}}  
-
\MBF{V_{i-1}^T}
\left[
  \begin{array}{c}
    \MBF{0} \\
    \MBF{\tilde{F}_{i-1}} \\ 
    \MBF{0} \\
  \end{array}
\right]
&=
\MBF{Y_{it}} 
+
\MBF{V_{i-1}^T}
\left[
  \begin{array}{c}
    \MBF{Y_{i-1t}} \\ 
    \MBF{0}\\
    \MBF{0}
  \end{array}
\right]
+
\MBF{V_{i-1}^T}
  \left[
    \begin{array}{c}
    \MBF{0}\\
    \MBF{0}\\
    \MBF{Y_{i-1b}}
  \end{array}
\right]
\\ &= 
\MBF{Y_{it}}
+
\MBF{V_{i-1t}^T}
\MBF{Y_{i-1t}}
+
\MBF{V_{i-1b}^T}
\MBF{Y_{i-1b}},
\end{split}
\end{gather}

\begin{gather} \label{gform2}
\begin{split}
i<p-1, \quad \MBF{G_{ib}} = 
\MBF{F_{ib}}  
-
\MBF{W_{i+1}^T}
\left[
  \begin{array}{c}
    \MBF{0} \\
    \MBF{\tilde{F}_{i+1}} \\ 
    \MBF{0} \\
  \end{array}
\right]
&=
\MBF{Y_{ib}} 
+
\MBF{W_{i+1}^T}
\left[
  \begin{array}{c}
    \MBF{Y_{i+1t}} \\ 
    \MBF{0}\\
    \MBF{0}
  \end{array}
\right]
+
\MBF{W_{i+1}^T}
\left[
  \begin{array}{c}
    \MBF{0}\\
    \MBF{0}\\
    \MBF{Y_{i+1b}}
  \end{array}
\right]
\\ &=
\MBF{Y_{ib} }
+
\MBF{W_{i+1t}^T}
\MBF{Y_{i+1t}  }
+
\MBF{W_{i+1b}^T}
\MBF{Y_{i+1b}},
\end{split}
\end{gather}

%%%%%%%%%%%%%%%%%%%%%%%%%%%%%%%%%%%%%%%%%%%%%%%%%%%%%%%%%%%%%%%%%%%%%%%%%%%%%%%%%%%%%%%%%%%%%%%%%%%%%%%%%
%%%%%%%%%%%%%%%%%%%%%%%%%%%%%%%%%%%%%%%%%%%%%%%%%%%%%%%%%%%%%%%%%%%%%%%%%%%%%%%%%%%%%%%%%%%%%%%%%%%%%%%%%%%%%%%%%%5
  At this point it should be noted that once the factorization stage done using our SPIKE implementation,
  the middle values of $\MBF{V_i^T}$ and $\MBF{W_i^T}$ are not available anymore, and they cannot then be used to construct  the components of
  $\MBF{G_i}$ in (\ref{gform1}) and (\ref{gform2}).
  Using the expression (\ref{eq:spikes}) for the spikes $\MBF{V_i}$ and $\MBF{W_i}$, $\MBF{G_i}$ can be rewritten as:
  
\begin{gather} \label{G_in_transpose_1}
i>1, \quad \MBF{G_{it}} 
= 
\MBF{F_{it}} - (\MBF{A^{-1}_{i-1}}
    \MBF{C_{i-1}}
	)\MBF{^T}
\left[
  \begin{array}{c}
    \MBF{0} \\
    \MBF{\tilde{F}_{i-1}} \\ 
    \MBF{0} \\
  \end{array}
\right]
= 
\MBF{F_{it}} - 
\left[
  \begin{array}{ccc}
    \MBF{\hat{C}_{i-1}^T} & 0 & \dots\\
    %\MBF{0} \\
  \end{array}
\right]
%\MBF{^T}
\MBF{A^{-T}_{i-1}}
\left[
  \begin{array}{c}
    \MBF{0} \\
    \MBF{\tilde{F}_{i-1}} \\ 
    \MBF{0} \\
  \end{array}
\right],
\end{gather}

\begin{gather}\label{G_in_transpose_2}
i<p-1, \quad \MBF{G_{ib}} 
= 
\MBF{F_{ib}} - (\MBF{A^{-1}_{i+1}}
    \MBF{B_{i+1}} 
	)\MBF{^T}
\left[
  \begin{array}{c}
    \MBF{0} \\
    \MBF{\tilde{F}_{i+1}} \\ 
    \MBF{0} \\
  \end{array}
\right]
= 
\MBF{F_{ib}} - 
\left[
  \begin{array}{ccc}
    \dots & \MBF{0} & %\\
    \MBF{\hat{B}_{i+1}^T} 
  \end{array}
	\right]%\MBF{^T}
\MBF{A^{-T}_{i+1}}
\left[
  \begin{array}{c}
    \MBF{0} \\
    \MBF{\tilde{F}_{i+1}} \\ 
    \MBF{0} \\
  \end{array}
\right].
\end{gather}

Overall, this approach is preferable to using the $\MBF{V_i}$ and $\MBF{W_i}$ matrices for two reasons.

  First,  as it can
be seen in Figure~\ref{reduced_transpose_1},  the top tip of $Y_1$ and the bottom tip of $Y_p$ 
make it through this transpose S-stage unchanged (resp. $Y_{1t}=F_{1t}$ and $Y_{pb}=F_{pb}$).
Therefore, the spikes $V_1$ and $W_p$ do not need to be formed during the factorization stage leading
to the load balancing optimization presented in Section~\ref{Load_balancinG_scheme_section} (i.e. the first and last partition can be chosen
bigger in size).

Second, $\MBF{G_{i+1t}}$ and $\MBF{G_{i-1b}}$ both require the same solve operation over the modified $\MBF{F_i}$ vectors,
\begin{equation}\label{eq:ATF}
\MBF{A_{i}^{-T}}
\left[
  \begin{array}{c}
    \MBF{0} \\
    \MBF{\tilde{F}_{i}} \\ 
    \MBF{0} \\
  \end{array}
\right].
\end{equation}
  Therefore, creating the $\MBF{G}$ vector in this manner incurs the cost of one large solve operation and two small multiplications per partition
 (since $\MBF{B_{i+1}}$ and   $\MBF{C_{i-1}}$ are mostly comprised of zeroes).
  This is likely to be less expensive than the cost of performing two large multiplications (if $V_i$ and $W_i$ were available).

Once the reduced system and $\MBF{G}$ vector have been constructed, all that remains in the S stage is to solve it.
Notably, this reduced system matrix is simply the transpose of the reduced system matrix used in non-transpose SPIKE given in (\ref{RED_SYS_RECA}) for four partitions.
In Section~\ref{RRS} a recursive method for solving the transpose reduced system will be presented.

%%%%%%%%%%%%%%%%%%%%%%%%%%%%%%%%%%%%%%%%%%%%%%%%%%%%%%%%%

\subsection{Transpose D stage} \label{DStage}

\begin{figure}[htbp] 
\begin{centering}
\includegraphics[keepaspectratio]{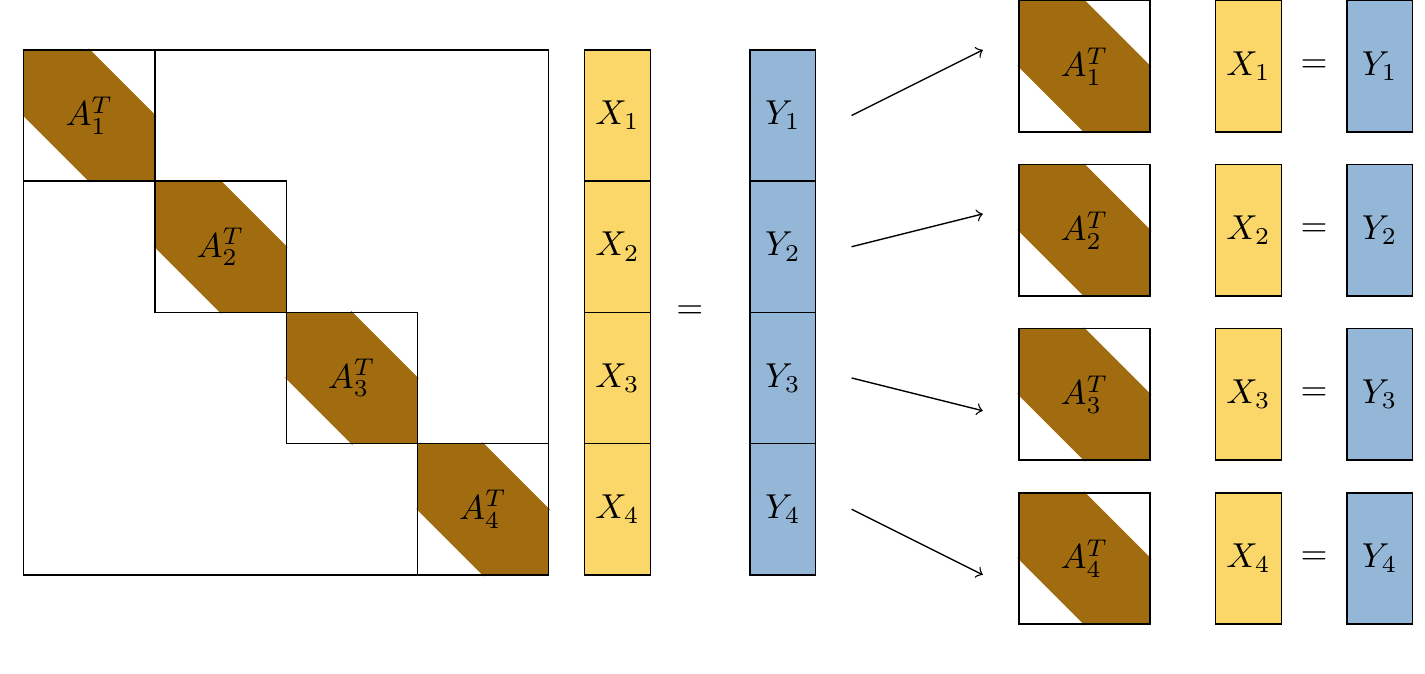} 
\caption{Transpose D stage}
\label{D_stage_split}
\end{centering}
\end{figure}

Because the partitions of the $\MBF{D}$ matrix are completely decoupled, performing this stage is much simpler than the S stage as illustrated in Figure~\ref{D_stage_split}.
The overall goal is to obtain $\MBF{X}$ in $\MBF{D^T}\MBF{X}=\MBF{Y}$.
In the S stage, it was shown that $\MBF{\tilde{Y}_i}=\MBF{\tilde{F}_i}$. Therefore, once the solutions of the reduced system $Y_{it}$ and $Y_{ib}$ are known, the
whole solution $X_i$ is simply retrieved as follows:

\begin{gather} \label{sweeps_in_D}
  \MBF{X_i} = 
  \MBF{A_i^{-T}}
  \left[
    \begin{BMAT}(b,14pt,14pt){c}{ccc}
      \MBF{Y_{it}} \\
      \MBF{\tilde{Y}_i} \\
      \MBF{Y_{ib}} 
    \end{BMAT}
  \right]
  =
  \MBF{A_i^{-T}}
  \left[
    \begin{BMAT}(b,14pt,14pt){c}{ccc}
      \MBF{Y_{it}} \\
      \MBF{\tilde{F}_i} \\
      \MBF{Y_{ib}} 
    \end{BMAT}
  \right].
\end{gather}

This concludes the description of the basic transpose SPIKE solver. 

  Similarly to the non-transpose case, optimizations are possible for transpose SPIKE to achieve the same computational costs
  reported in Table~\ref{Computation_Cost_Table}
  for the total number of solve sweeps depending of the type of partition \cite{master}.
  The fact that it is not necessary to generate the full $W$ spike for the first partition and $V$ spike for the last partition,
   allows for the creation of a $2\times2$ ``transpose'' kernel, which can be used for developing
   a flexible threading strategy applied to transpose SPIKE similar to the one presented in Section~\ref{Increased_parallelism_for_Recursive_Spike}.

\subsection{Transpose recursive reduced system} \label{RRS}

In Section~\ref{recursive_reduced_section}, a description of the recursive method of solving the reduced system was described. 
Because the reduced system of transpose SPIKE is simply the transpose of the original reduced system,
it suffers from the same problem: increasing the number of partitions increases the size of the reduced system. 
Therefore, a recursive method for solving the reduced system is also required for the transpose case.

For the transpose reduced system, we aim at reusing the recursive factorization performed for the non-transpose case.
The result from a second level of SPIKE DS factorization applied to the original
  reduced system was given in (\ref{RED_SYS_RECB}) (using half the number of partitions):

\begin{equation}
\label{transred1}
\MBF{S^{[1]}} = \MBF{D^{[1]}} \MBF{S^{[2]}},
\end{equation}
and this process can be repearted on the new generated spike matrix until only two partitions are left, i.e.
\begin{equation}
\label{transred2}
\MBF{S}^{[i]} = \MBF{D^{[i]}} \MBF{S^{[i+1]}}.
\end{equation} 

\begin{figure}[htbp] 
\begin{centering}
\includegraphics[keepaspectratio,height=\textheight]{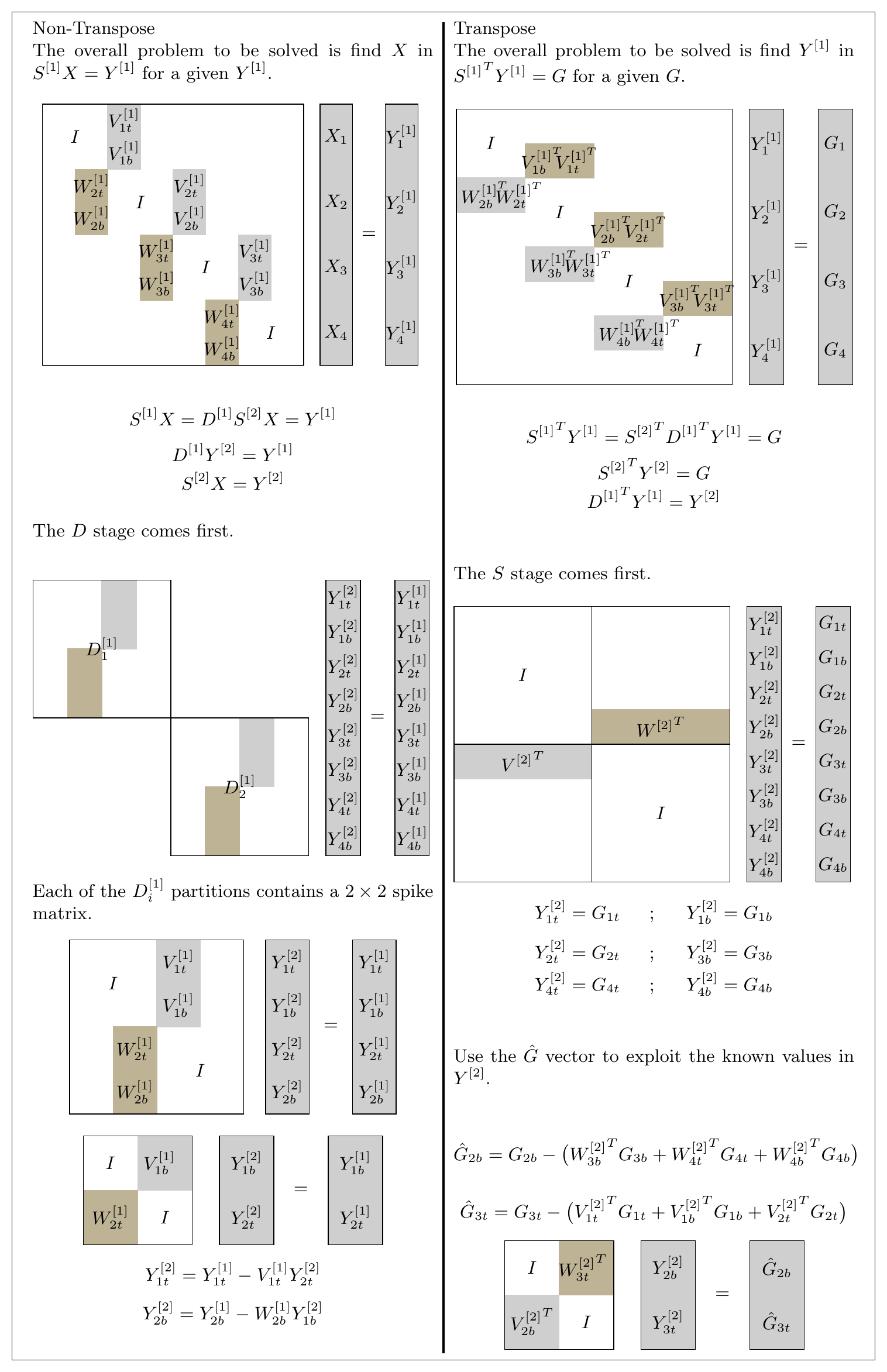} 
\caption{SPIKE four partition recursive reduced system solve, part 1}
\label{reduced_example_solve_A}
\end{centering}
\end{figure}

\begin{figure}[htbp] 
\begin{centering}
\includegraphics[keepaspectratio,height=\textheight]{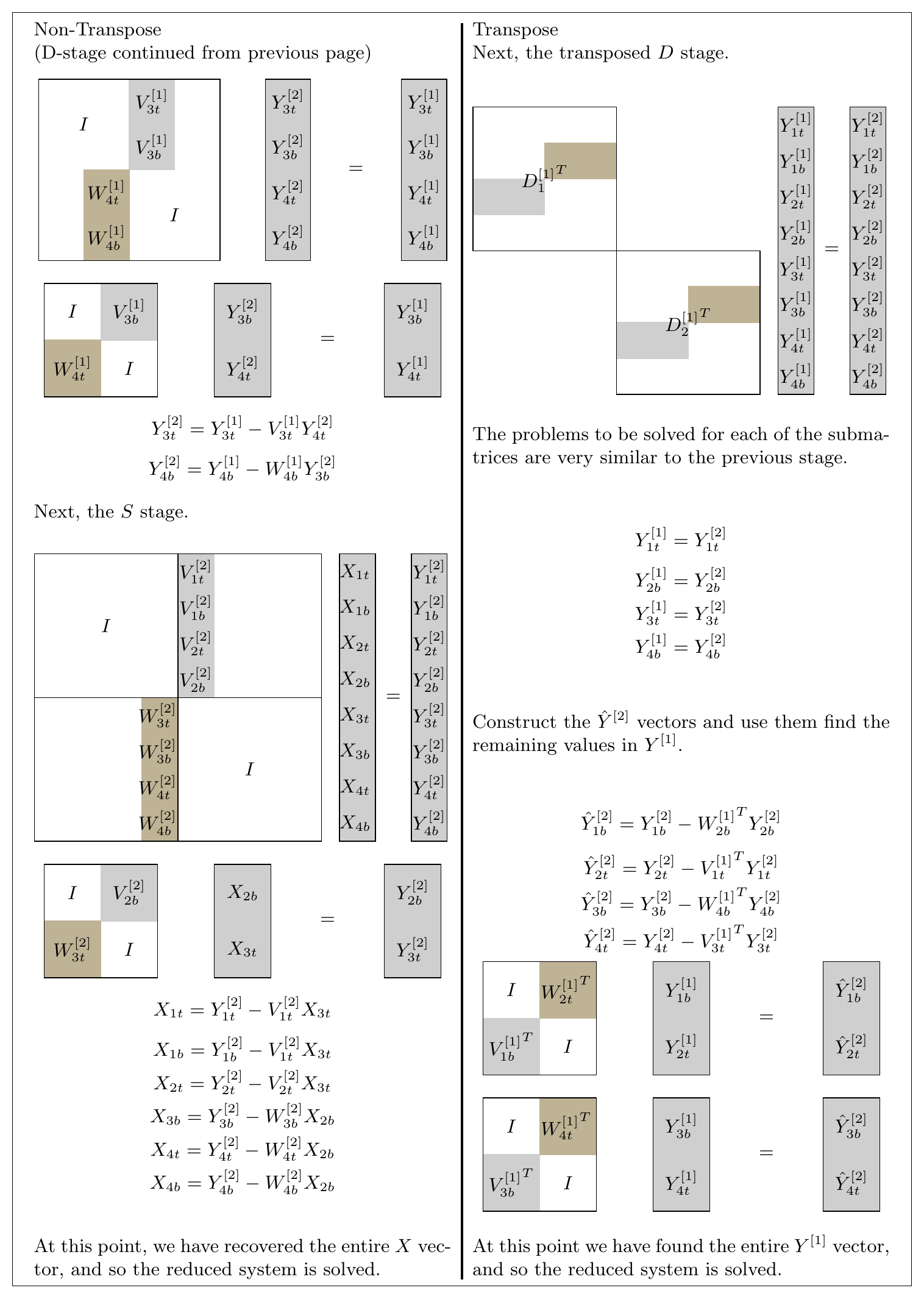} 
\caption{SPIKE four partition recursive reduced system solve, part 2}
\label{reduced_example_solve_B}
\end{centering}
\end{figure}

With each step of this recursion, the number of partition is divided by two and the size of the partitions doubles. 
  If $p$ is the number of partitions into which the original matrix was broken, the process can be repeated in $\RECLEVELS=\log_2(p)$ times
  \cite{Polizzi:2006}. It comes:
\begin{gather}
\MBF{{S}^{[1]}}  = \left(\Pi_{i=1}^{\RECLEVELS-1} \MBF{D^{[i]}}\right) \MBF{S^{[\RECLEVELS]}},
\end{gather} 
where $\MBF{S^{[\RECLEVELS]}}$ has only two partitions left.
For the transpose case, we have $\MBF{{S}^T}\MBF{Y^{[1]}}=\MBF{G}$ (see Fig. \ref{reduced_transpose_1}),
so we may perform the transpose operation on the series of products above:
\begin{gather}
\MBF{{S}^T}  = \MBF{S^{[\RECLEVELS]}}^T\left(\Pi_{i=\RECLEVELS-1}^{1} \MBF{D^{[i]}}^T\right).
\end{gather} 
This could be thought of as performing the original, non transpose, reduced system solve, but with the solve stages in reverse. 
The operation to be performed is:
\begin{gather}
\MBF{Y^{[1]}} = \MBF{{S}^{-T}} \MBF{G} = \left(\Pi_{i=1}^{\RECLEVELS-1} \MBF{D^{[i]}}^{-T}\right)\MBF{S^{[\RECLEVELS]}}^{-T}\MBF{G}.
\end{gather}
  The full process of solving the reduced system using four partitions, is shown in
Figures~\ref{reduced_example_solve_A} and~\ref{reduced_example_solve_B} where non-transpose and transpose cases are detailed side-by-side.

%%%%%%%%%%%%%%%%%%%%%%%%%%%%%%%%%%%%%%%%%%%%%%%%%%%%%%%%%%%%%%%%%%%%%%%%%%%%%%%%%%%%%%%%%%%%%%%%%%%%%%%%%%%

%%%%%%%%%%%%%%%%%%%%%%%%%%%%%%%%%%%%%%%%%%%%%%%%%%%%%%%%%%%%%%%%%%%%%%%%%%%%%%%%%%%%%%%%%%%%%%%%%%%%%%%%%%%%%%%%%%%%%%%%%%%%

\subsection{Transpose solver performance}

\begin{figure}[htbp] 
\includegraphics[keepaspectratio,width=\textwidth]{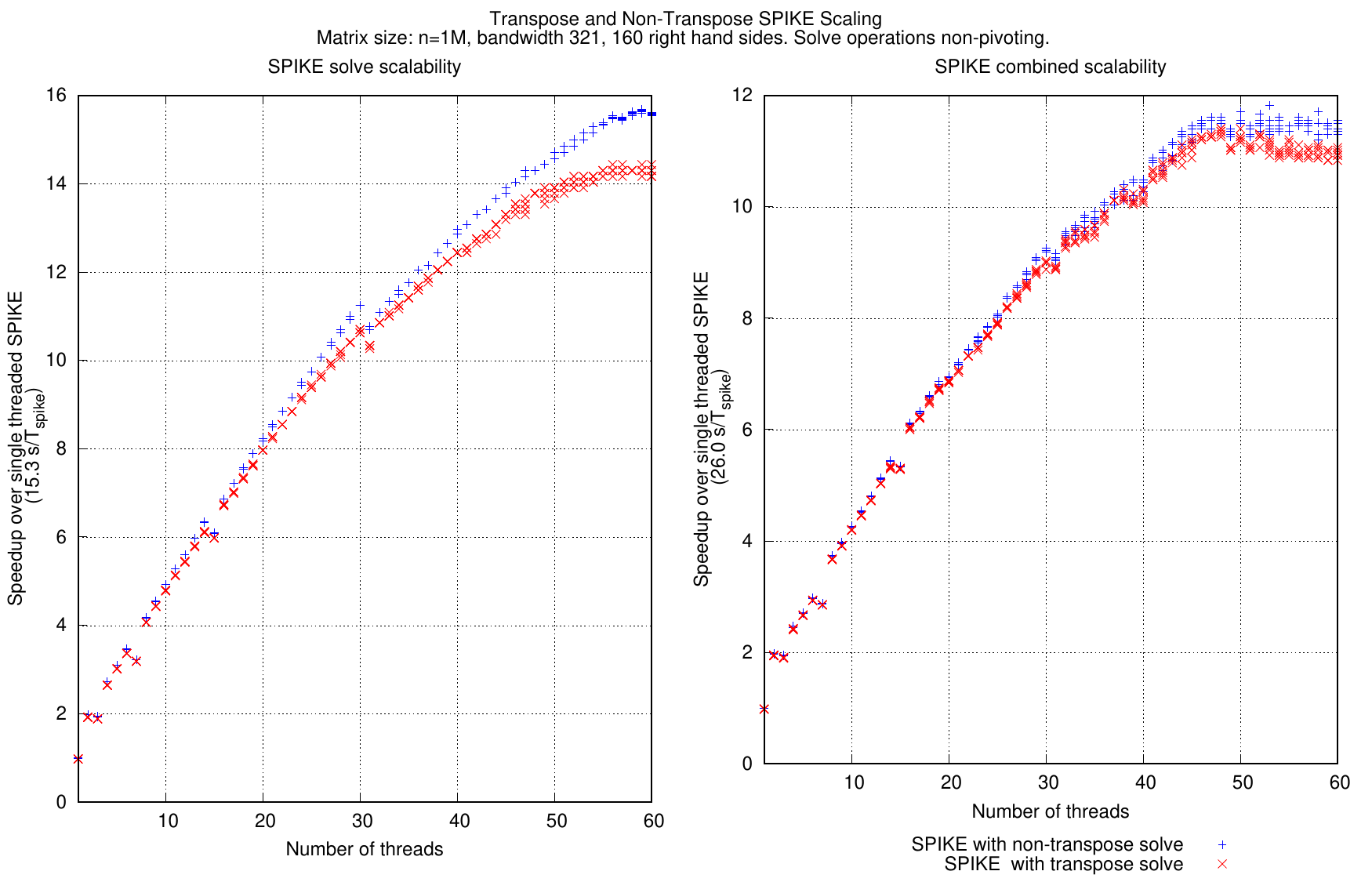} 
\caption{Computation time comparisons}
\label{scalabilitytrans160}
\end{figure}

Figure~\ref{scalabilitytrans160} shows the solve stage, as well as overall, scaling compared to the single-threaded non-pivoting non-transpose solver. 
This base solver was chosen to make a one-to-one comparison with the non-transpose solver.
Because the factorization is reused for both the transpose and non-transpose problem, factorization time is not shown. 

The transpose option has little effect on performance. 
There is a very slight performance loss in the overall case, and a more noticeable one when just looking at the solve stage. 
However, in either case, the loss of performance generally occurs well past the point where diminishing returns have already set in, and does not appear to degrade overall performance significantly. 

%%%%%%%%%%%%%%%%%%%%%%%%%%%%%%%%%%%%%%%%%%%%%%%%%%%%%%%%%%%%%%%%%%%%%

\section{An efficient pivoting scheme} \label{Pivoting_Spike_Section} 

The standard \lapack libraries use partial pivoting to increase the numerical stability of the solve operation \cite{HighamLUAccuracy}. 
Partial pivoting operates by exchanging rows when the pivot element is selected, placing the greatest element in the column on the diagonal.
This decreases the loss of accuracy caused by rounding, and reduces the chances of selecting zero as the pivot element.

As originally described in \cite{Polizzi:2006}, the recursive SPIKE algorithm is
  using non-pivoting factorization schemes along with a diagonal boosting strategy.
With diagonal boosting, a small value is added to zero-pivots when they are discovered, resulting in an approximate factorization.
SPIKE would then operate as a good preconditioner since few iterative refinements are generally needed to reach convergence.
Interestingly, the diagonal boosting strategy could also be a viable option
in the case where partial pivoting fails (since full-pivoting solver are not readily available).
The non-pivoting option in SPIKE helps maintaining the banded structure of the matrix, which simplifies the implementation of the algorithm
and improves performance of the factorization stage. Although, in the large majority of cases
  zero-pivot are rare in double precision arithmetic (so boosting may not occurred),
  partial pivoting for SPIKE may become a necessity if the matrices are not very well conditioned.
In addition, an efficient  partial pivoting SPIKE solver could allow better one to one comparison with LAPACK LU solver.

\subsection{Pivoting LU factorization}

\begin{figure}[htbp] 
\begin{centering}
\includegraphics[keepaspectratio,width=0.5\textwidth]{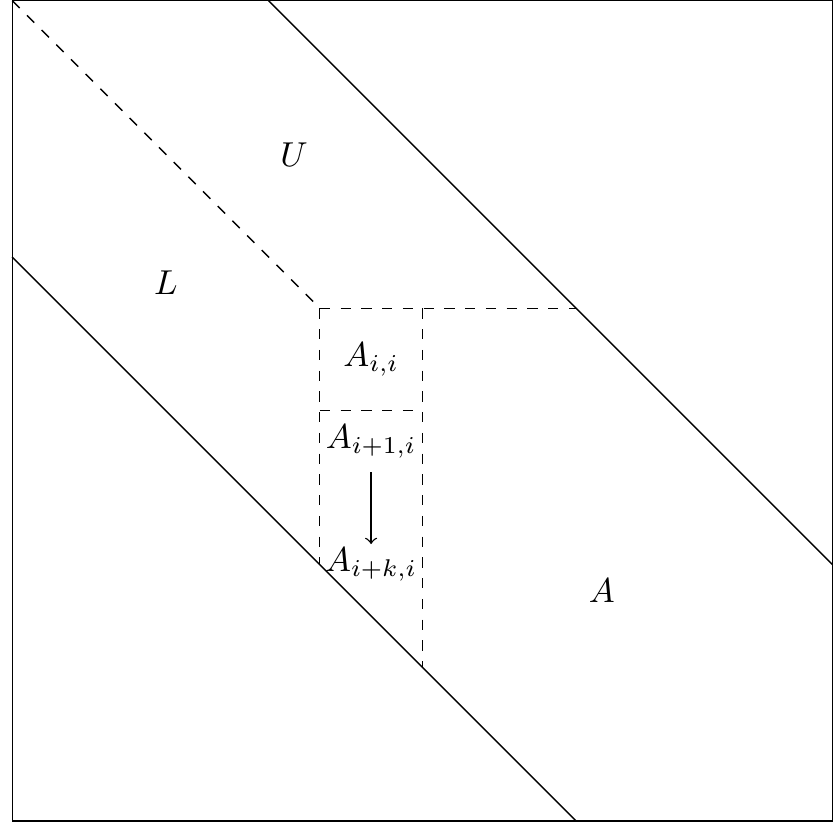} 
\caption{LU factorization in intermediate step $i$, potential pivot elements are $\MBF{A_{i,i}}$ to $\MBF{A_{i+k,i}}$.}
\label{pivot_cartoon}
\end{centering}
\end{figure}

The algorithm implemented for the \lapack LU factorization  is essentially similar to the Doolittle algorithm. 
In particular, the L and U matrices are crafted column-by-column, progressing from left to right along the diagonal. 
As a result, the only legitimate selections for pivot rows are those below the diagonal as shown in Figure~\ref{pivot_cartoon}. %\cite{HighamLUAccuracy}.
In addition, the row selected must have a non-zero value, restricting the choices to those within the band. 
So, the only possible candidates for row exchange are those rows between the diagonal and the bottom most subdiagonal element in the current column, which means that a given permutation, at most, moves a row $k$ places up.

Because partial pivoting is implemented as a series of row exchange permutations, it can be represented as left-multiplication of a permutation matrix, $\MBF{P}$. 
Actually, the permutations are implemented inside of the L-sweep. 
As a result, the pivoting LU factorization and solve operations can be represented as follows, for some arbitrary collections of vectors $\MBF{F}$ and $\MBF{X}$:

\begin{equation}
\MBF{PA} = \MBF{LU} 
\end{equation}
\begin{equation}
\MBF{A^{-1}}\MBF{F} 
=\MBF{U^{-1}L^{-1}P}\MBF{F}
=\MBF{U^{-1}}\MBF{\left(L^{-1}PF \right)}
=\MBF{X}   
\end{equation}

In other words, the effect of partial pivoting on the collection of vectors $\MBF{F}$ is the same as the effect on the matrix $\MBF{A}$. 
An element in $\MBF{F}$ may be moved at most $k$ places upwards.  In the context of SPIKE,
this will allow us to continue using the optimizations which exploit the triangular shape of the $\MBF{L}$ and $\MBF{U}$ matrices.  
These optimizations are described in Section~\ref{opt_sweeps}, and the related operations are performed for equations~(\ref{V_formula_end}) and~(\ref{one}).
First, looking at~(\ref{V_formula_end}), the original equation was

\begin{equation} 
\MBF{V_i} = \MBF{A_i^{-1}} \MBF{B_i}
=
\MBF{U_i^{-1}}\MBF{L_i^{-1}} 
\left[\begin{BMAT}(b,18pt,18pt){cc}{cc}
 \MBF{0}  &\MBF{0} \\
 \MBF{\hat{B}_i}  &\MBF{0} 
\end{BMAT}\right].
\end{equation}
The permutation matrix must now be inserted as follows  
\begin{equation} 
\MBF{V_i} = \MBF{A_i^{-1}} \MBF{B_i}
=
\MBF{U_i^{-1}}\MBF{L_i^{-1}}\MBF{P_i}
\left[\begin{BMAT}(b,18pt,18pt){cc}{cc}
 \MBF{0}  &\MBF{0} \\
 \MBF{\hat{B}_i}  &\MBF{0} 
\end{BMAT}\right].
\end{equation}
When performing solve operation with $\MBF{L}$, we may simply break up the zero-matrices as follows:
\begin{equation} 
\MBF{L_i^{-1}}\MBF{P_i}
\left[\begin{BMAT}(b,18pt,18pt){cc}{cc}
 \MBF{0}  &\MBF{0} \\
 \MBF{\hat{B}_i}  &\MBF{0} 
\end{BMAT}\right]
=
\MBF{L_i^{-1}}\MBF{P_i}
\left[\begin{BMAT}(b,18pt,18pt){cc}{ccc}
 \MBF{0}  &\MBF{0} \\
 \MBF{\hat{0}}  &\MBF{\hat{0}} \\
 \MBF{\hat{B}_i}  & \MBF{0}
\end{BMAT}\right],
\end{equation}
where $\MBF{\hat{0}}$ is a matrix with $k$ rows. 
Now, we may begin the L-sweep at the top of $\MBF{\hat{0}}$, and any pivoted rows of $\MBF{B}$ will still be involved in the solve operation. 
From here, the operations may continue as in non-pivoting SPIKE. 

\subsection{Pivoting UL factorization}
There is no UL factorization specified in \lapack.  
However, a efficient UL factorization and solve is necessary to reduce the number of solve sweeps used in the last SPIKE partition, as shown in Section~\ref{opt_sweeps}.
Specifically, we require the ability to obtain the topmost elements of $\MBF{W_p}$ without using any large sweeps, and limit the contamination caused by the reduced system to the topmost elements of $\MBF{Y_p}$. 

Implementing a pivoting UL factorization with performance comparable to, for example, Intel MKL is clearly beyond the scope of this project. 
Instead we use a permutation to effectively obtain a UL factorization using the native \lapack LU factorization.
The permutation matrix, given as Q below, has ones on the anti-diagonal. 
 
\begin{equation} 
\MBF{Q}=
\left[\begin{BMAT}(b,18pt,18pt){ccccc}{ccccc} 
   &        &   &        & 1 \\
   &        &   & \udots &   \\
   &        & 1 &        &   \\
   & \udots &   &        &   \\
 1 &        &   &        &   
\end{BMAT}\right]
\label{Qmat}
\end{equation}

$\MBF{Q}$ has the property that pre-multiplying some matrix by $\MBF{Q}$ reverses the order of the rows of that matrix, and post-multiplying a matrix by $\MBF{Q}$ reverses the order of the columns. 
It is also orthogonal and symmetric; $\MBF{Q}=\MBF{Q^T}=\MBF{Q^{-1}}$. 
Thus, a given matrix solve problem may be rewritten as follows

\begin{equation}
\MBF{AX}=\MBF{F}=\MBF{QQAQQX}=\MBF{Q(QAQ)QX}
\end{equation}

it comes
\begin{equation}
  X=Q(QAQ)^{-1}QF
\end{equation}

Because both the rows and columns of $\MBF{QAQ}$ have been reversed, this matrix is still banded.
So, it still may be operated upon using the standard pivoting LU factorization. In addition, 
 the topmost elements of $\MBF{F}$ becomes the bottom most elements of $\MBF{QF}$.
As a result, the successive permutations and triangular solves
  can be performed from right to left, as follows:

\begin{equation}
X=Q\MBF{\left(U^{-1}\Big(L^{-1}\big(P(QF)\big)\Big)\right)}.
\end{equation}

Thus, the structure of the collections of vectors used for the final partition is essentially the same as that of the vectors used in the first partition. 
$\MBF{QW_p}$ has the same essential shape as $\MBF{V_1}$.
And so, we may reuse the same optimizations for the final partition as were used for the first. 

Finally, it is possible to perform the pivoting UL in place using the pivoting LU factorization,
  by explicitly moving the elements of the matrix and vectors around in memory. 
The computational and memory cost of this reordering is significantly less than that of the factorization of the full permuted matrix $QAQ$.  
 A dedicated pivoting UL factorization would be the best alternative since our current approach for UL factorization
 could impact scalability noticeably (as it will be shown in benchmarking). However, this method does not prevent progress completely.

\subsection{Performance measurements}
\subsubsection{Computation Time} 
The purpose of pivoting SPIKE is to reduce the accuracy loss associated with using a non-pivoting solver,
while retaining some of the performance advantage over a pivoting one. 
So, the relevant metrics are the computation time, scaling, and the residual produced. 
The use of a pivoting solver has two noticeable performance impacts.
First, during the 
factorization, the pivot element is selected by scanning through the 
column and locating the element with the greatest magnitude. This 
scanning process occurs independent of the  diagonal dominance.
Second, when the matrix is not diagonally dominant, there is a cost associated with performing the pivoting operation. 

For the sake of these comparisons, it is useful to vary both the number of threads and the diagonal dominance of the matrix. 
As a slight extension to the concept of a diagonally dominant matrix, let us define $DD$, the 'degree of diagonal dominance,' as the following:
\begin{equation}
DD=\min_{i \in 1\dots n}{\left(\frac{A_{ii}}{\sum_{j \neq i}A_{ji}}\right)}
\end{equation}

A diagonally dominant matrix would have $DD \geq 1$.
To generate matrices with a desired value for $DD$, the following procedure has been used: 
Each element within the non-zero band of the matrix has been filled with random values using the \lapack DLARNV command. 
Then, the columns are summed and multiplied by the desired value for $DD$ and the result is placed on the diagonal. 

\begin{figure}[htbp] 
\includegraphics[keepaspectratio,width=\textwidth]{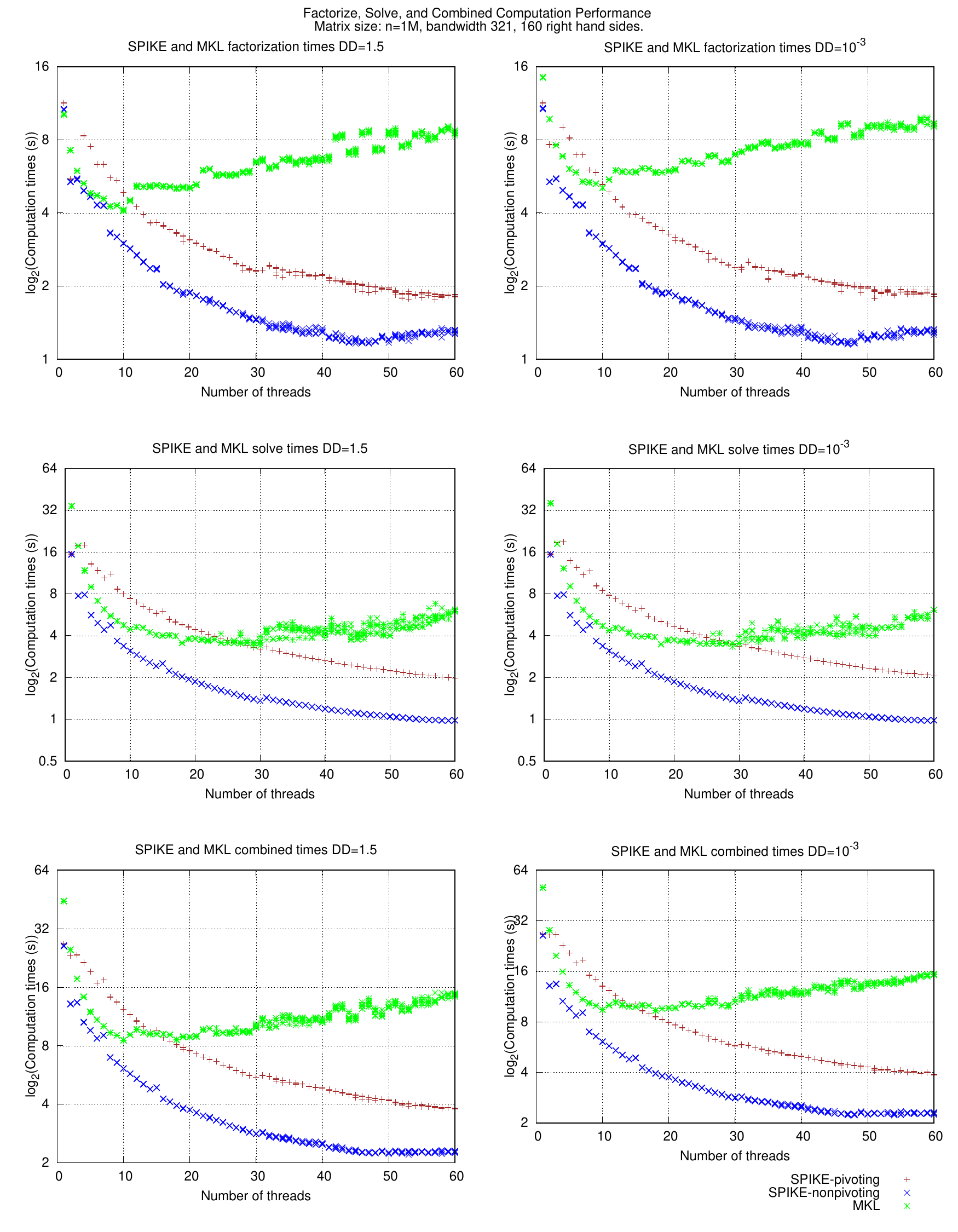} 
\caption{Computation time comparisons}
\label{scalabilitypivot160}
\end{figure}

Figure~\ref{scalabilitypivot160} shows the overall performance comparisons for non-pivoting SPIKE, pivoting SPIKE, and MKL.
Note that computation time is plotted on a log scale to retain the visibility of performance changes for large numbers of threads.
The hardware and software used for these runs were detailed in Section~\ref{sec:machine}.
Two matrix configurations are used, one in which the matrix is diagonally dominant ($DD=1.5$), and one in which it is not ($DD=10^{-3}$).
Non-pivoting SPIKE  clearly demonstrates the best performance.
Pivoting SPIKE and MKL perform well in different conditions, with MKL obtaining a noticeable advantage for low numbers of threads
-- the additional cost of not having a dedicated and optimal pivoting UL factorization is a likely cause of this issue (involving also an additional permutation
in the solve stage).
SPIKE improves in performance as the number of threads increases. 
In particular, the MKL factorization stage does not scale well beyond 10 threads on this machine, likely because at this point the computation
begins to access additional processor packages.
Overall, it would appear that the SPIKE decomposition technique is quite helpful in improving performance scalability.

\subsubsection{Precision}
\begin{figure}[htbp] 
\includegraphics[keepaspectratio,width=\textwidth]{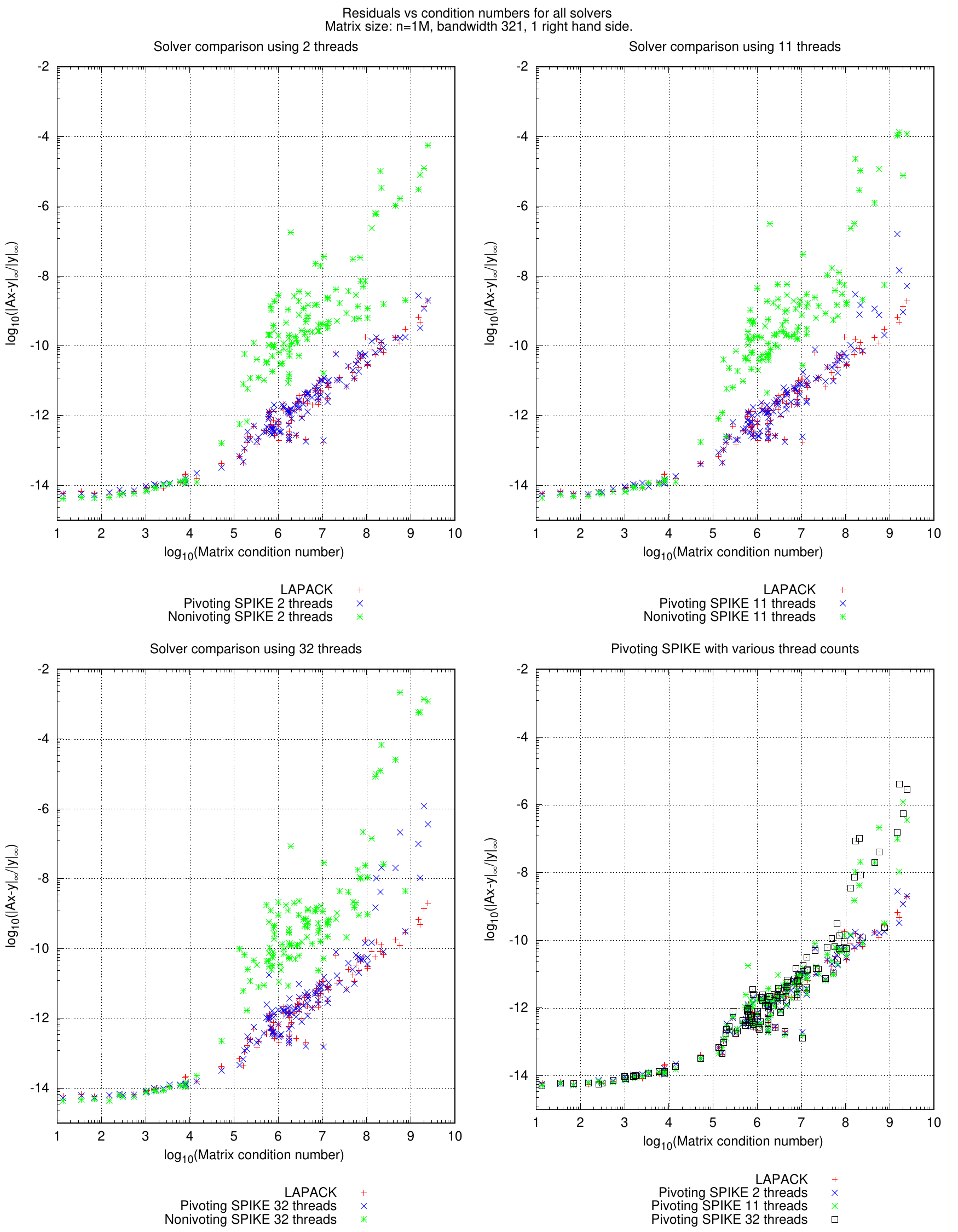} 
\caption{Condition number and residual relationship}
\label{residuals}
\end{figure}

Figure~\ref{residuals} shows the numerical accuracy advantages of pivoting SPIKE, by comparing the residual produced to the condition number. 
Matrices are produced in the same manner as the preceding section, and condition number of estimated by the \lapack function DGBCON. 
All computations are performed in double precision.

The top-left, top-right, and bottom-left quadrants of the figure compare the three solvers.
In the top-left quadrant it can be seen that, with two-partitions, pivoting SPIKE produces residuals indistinguishable from \lapack.
Results for non-pivoting SPIKE are also comparable for condition numbers less than $10^5$. The residuals start increasing after this point for all solvers,
with a noticeable much higher increase for non-pivoting SPIKE. 
In the top-right and bottom left quadrants we see some loss of accuracy for the pivoting SPIKE, particularly as the condition number becomes very large. 

The bottom-right quadrant shows a comparison of pivoting solvers for all thread counts used. 
Viewing this chart, it becomes apparent that there are three relevant ranges for the computation. 
For condition numbers in the range of $1$ to $10^{5}$, all of the solvers perform well. 
For condition number in the $10^{5}$ to $10^{8}$, the residuals produced by the pivoting solvers are essentially identical. 
Finally, for condition numbers greater than $10^{8}$ there is some loss of precision for pivoting SPIKE based on the number of partitions used. 

In summary, the residuals produced by the pivoting SPIKE solver are a significant improvement over non-pivoting SPIKE for poorly conditioned matrix.
There are cases where the pivoting SPIKE solver loses precision as the number of partitions increases, but for the range
of interesting problems the pivoting SPIKE solver precision is indistinguishable to the pivoting \lapack solver. 

\section{Conclusion}

A feature complete recursive SPIKE algorithm has been presented. Three enhancements for SPIKE have been shown,
achieving near feature-parity with the standard \lapack banded matrix solver.
In particular, both the transpose solve option and the partial pivoting option, provide standard capabilities
found in \lapack solvers.
Transpose solve operation allows improved algorithmic flexibility and efficiency by eliminating the need for
an additional transpose factorization. Pivoting operation provides a convenient middle-ground between the numerical accuracy of the standard LAPACK
solver and the extreme scalability of the standard SPIKE algorithm. 

All algorithms have been implemented with a flexible threading scheme that allows the effective utilization of %nearly
any number of threads, overcoming a previous known limitation of the recursive SPIKE scheme.
In addition, the per-partition performance has been characterized, resulting in a simple load balancing equation
controlled by a single machine specific parameter.
With the addition of these features and demonstrated performance advantages,
it is our hope that the new SPIKE-OpenMP library \cite{spike} may be considered a drop-in replacement
for the standard \lapack banded factorize and solve operations.

\begin{acks}
   This work was supported by National Science Foundation grant CCF-\#1510010.
\end{acks}

\bibliographystyle{ACM-Reference-Format-Journals}
\bibliography{citations}
                             % Sample .bib file with references that match those in
                             % the 'Specifications Document (V1.5)' as well containing
                             % 'legacy' bibs and bibs with 'alternate codings'.
                             % Gerry Murray - March 2012

% History dates
%\received{February 2007}{March 2009}{June 2009}

% Electronic Appendix
%\elecappendix

%\medskip

\end{document}